\documentclass[12pt,final]{amsart}

\usepackage[english]{babel}
\usepackage{graphicx}

\usepackage{amsmath,amsthm,amsfonts,amssymb}
\usepackage{csquotes}
\usepackage[round]{natbib}
\usepackage{multirow}
\usepackage{xcolor}
\usepackage{url}
\usepackage{algorithm2e}

\newcommand\R{\mathbb{R}}

\newcommand\Ld{\mathbb{L}^2([0,1])}

\title[RSM for functional data]{Local optimization of black-box function with high or infinite-dimensional inputs. \\Application to nuclear safety.}
\author[A. Roche]{Angelina Roche\\}
\address{MAP5 UMR CNRS 8145, Universit\'e Paris Descartes and CEREMADE UMR CNRS 7534, Universit\'e Paris Dauphine.}
\email{angelina.roche@dauphine.fr}
\urladdr{https://www.ceremade.dauphine.fr/~roche/}
\date{\today}
\begin{document}
\begin{abstract} An adaptation of Response Surface Methodology (RSM) when the covariate is of high or infinite dimensional is proposed, providing a tool for black-box optimization in this context. We combine dimension reduction techniques with classical multivariate Design of Experiments (DoE). We propose a method to generate experimental designs and extend usual properties (orthogonality, rotatability,...) of multivariate designs to general high or infinite dimensional contexts.  Different dimension reduction basis are considered (including data-driven basis). The methodology is illustrated on simulated functional data and we discuss the choice of the different parameters, in particular the dimension of the approximation space. The method is finally applied to a problem of nuclear safety.  

\end{abstract}
\maketitle
$\;$


$\;$

\section{Introduction}
Black-box optimization problems arises in many applications, for instance when one wants to optimise an output of a computer code or in real-life experiments such as crash test, chemical reactions,  medical experiments... In more and more applications, the input is high-dimensional, or even infinite-dimensional (time or space dependent).  In this paper, our aim is to minimise an unknown function $m:\mathcal R\to \mathbb R$ where $\mathcal R$ is a subset of a separable Hilbert space $(\mathbb H,\langle\cdot,\cdot\rangle,\|\cdot\|)$, which can be e.g. $\mathbb R^d$ or a function space. The function $m$ is unknown, but noisy evaluations of $m$ are available. We suppose that each evaluation of $m$ is costly, thus the aim is to be as close as possible to an optimum with a given (low) number of evaluations of $m$.

  In this context, surrogate-based approaches, such as those based on response-surface methodology, kriging, radial basis functions, splines or neural networks are commonly used \citep{queipo_surrogate_2005,simpson_metamodels_2001}. In Response Surface Methodology (RSM), the function $m$ is locally approximated by a polynomial regression model, typically with order 1 or 2  \citep[see e.g. the review of][]{khuri_response_2010}.  With the information given by the evaluation of the function $m$ on design points chosen by the user, least-squares estimates of the model parameters provide a local approximation of the surface $y=m(x)$. This local approximation can be used, for instance, to approximate the gradient or to locate a critical point of the surface. The efficiency of the method rests mainly on the choice of the design of experiments. Hence, a lot of research has been done in order to find sets of points giving the best possible precision of the fitted surface with few evaluations of $m$ on the design points. We refer to \citealt[pp.~131--133]{khuri_response_2010} for the description of the most common response surface designs. This subject is still an active research field (see e.g. \citealt{georgiou_class_2014}). However, these design generating methods are not tractable when the inputs are high-dimensional (for instance, $\mathbb H=\mathbb R^d$ with $d\geq 100$) and can not be directly defined with infinite-dimensional inputs.  
  
Projection-based dimension reduction techniques has become a powerful tool in high-dimensional statistics and are the main tool of most of the methods used to treat functional data. Usually, the data are projected into a subset $S=\text{Vect}\{\varphi_1,...,\varphi_d\}\subset\mathbb H$ of reduced dimension. In functional data analysis, $\mathbb H$ is a function space and fixed basis such as Fourier basis, spline or wavelet basis are commonly used, exploiting the regularity properties (smoothness for instance) of the functions in the sample. Another interesting approach consists in using the information given in a learning sample of pairs $\{(X_i,Y_i), i=1,...,n\}$ to generate the directions $\varphi_1,...,\varphi_d$. Among them the approaches based on the principal components  are the most common: PCA \citep{hall_principal_2011}, sparse PCA \citep{zou_sparse_2006,qi_sparse_2015}, regularized PCA \citep{rice_estimating_1991,lee_estimating_2002,ramsay_functional_2005}. Another approach to obtain interesting data-driven basis is Partial Leasts Squares regression \citep[PLS,][]{wold_soft_1975,preda_pls_2005}. The main advantage of PLS is that the directions $\varphi_1,...,\varphi_d$  are chosen in $\mathbb H$ so as to maximize the explained variance of $Y$. Hence, PLS uses the information of the whole learning sample $\{(X_i,Y_i), i=1,...,n\}$ whereas the PCA basis is generated only with the $X_i$'s.  

 The present paper proposes an adaptation of Response Surface Methodology when the input is in a general Hilbert space $\mathbb H$, via dimension reduction techniques. The specificities of the framework are explained in Section~\ref{sec:models}. In Section~\ref{sec:gen_func_DOE}, we provide a way to generate RSM design of experiments based on dimension reduction tools. A simulation study is presented in Section~\ref{sec:num_exp}. The method is finally applied in Section~\ref{sec:CEA} to a nuclear safety problem.
 
\section{High-dimensional and functional context}
\label{sec:models}
We suppose here that our real response $y$ depends on a variable $x$ in an infinite or high-dimensional space $\mathbb H$ equipped with a scalar product $\langle\cdot,\cdot\rangle$. For instance, $\mathbb H$ may be the space $\mathbb R^d$ equipped with its usual scalar product $ \langle x,y\rangle=\sum_{j=1}^d x_jy_j$ or a function space, such that $\mathbb L^2(I)$ with $I$ a measurable subset of $\mathbb R^d$ and $\langle f,g\rangle=\int_If(t)g(t)dt$. We denote by $\|\cdot\|$ the associated norm ($\|x\|^2=\langle x,x\rangle$ for all $x\in\mathbb H$).

We suppose here that the function $m$  is sufficiently smooth, so that the surface $y=m(x)$ can be approximated reasonably by a first or second-order surface. We consider here generalizations of the classical multivariate first and second-order models. Higher-order polynomial models, or even generalized linear models \citep{muller_generalized_2005,khuri_overview_2001}, while less standard, could also be considered similarly. 

\subsection{First-order model}
We define first-order models in the following form
 \[y: = \alpha+\langle\beta,x\rangle + \varepsilon,\]
 with $\alpha\in\R$, $\beta\in\mathbb H$ and $\varepsilon\sim\mathcal N(0,\sigma^2)$. If $\mathbb H$ is a function space, this model is known as \textit{functional linear model} \citep{ramsay_tools_1991,cardot_functional_1999} and has been widely studied (see \citealt{cardot_functional_2011} for a recent overview or \citealt{brunel_nonasymptotic_2016} for a recent work on this subject). Moreover, it is known that, if $\mathbb H$ is of high or infinite dimension, least-squares estimators of the slope parameter $\beta$ are, in general, unstable. Hence, precautions must be taken, either in the choice of the design (see Section~\ref{sec:gen_func_DOE}), or in the estimation method for the parameter $\beta$ which is an ill-posed inverse problem \citep{engl_regularization_1996}. The interest of this model is that, if $m$ is differentiable, for all $x_0\in\mathbb H$, $m(x)=m(x_0)+\langle m'(x_0),x-x_0\rangle+o(\|x-x_0\|)$, where $m'(x_0)$ is the gradient of $m$ at the point $x_0$, which implies that, if $x$ is sufficiently close to $x_0$, an estimator of $\beta$ will estimate the gradient of the surface $y=m(x)$ around the point $x_0$. 
 
 \subsection{Second-order model}
A second-order model can also be written as follows
 \begin{equation}\label{eq:functionalSOM}
 y: = \alpha+\langle\beta,x\rangle+\frac{1}{2}\langle H x,x\rangle + \varepsilon,\end{equation}
where $\alpha\in\R$, $\beta\in\mathbb H$ and $H:\mathbb H\to\mathbb H$ is a linear self-adjoint operator. Up to our knowledge, this model (at least in this particular form) has not been studied yet in the literature of functional data analysis. As for the first-order linear model, classical least-squares estimation is not a good choice and we have to be careful either in the choice of the design, or in the estimation method.   If $\mathbb H=\mathbb R^d$, the operator $H$ is a symmetric and positive semi-definite matrix $H=\left(h_{i,j}\right)_{1\leq i,j\leq d}$ and the model~\eqref{eq:functionalSOM} is the classical second-order multivariate linear model which can be written
\begin{equation*}
y=\alpha+\sum_{j=1}^d\beta_jx_j+\sum_{\substack{j,k=1\\j<k}}^dh_{j,k}x_jx_k+\sum_{j=1}^dh_{j,j}x_j^2+\varepsilon.
\end{equation*}

As for the second-order model, if $m$ is twice differentiable and if the design points are sufficiently close to $x_0$, an estimator of $\beta$ is an estimator of the gradient of $m$ and an estimator of $H$ gives an estimator of the Hessian matrix of $m$. In particular, if $x_0$ is close to a critical point, then the estimator of $\beta$ is close to 0 and the estimator of $H$ may help to precise the exact location and the nature of the critical point. 

 \section{Generation of Design of Experiments}
\label{sec:gen_func_DOE}

\subsection{General principle}

 The method is based on dimension reduction coupled with classical multivariate designs. The main idea is the following: suppose that we want to generate a design around $x_0\in\mathbb H$, we choose an orthonormal basis $(\varphi_j)_{j\geq 1}$ of $\mathbb H$, a dimension $d$ and a $d$-dimensional design $\{\mathbf x_i, i=1,\hdots,n\}=\{(x_{i,1},\hdots,x_{i,d}), i=1,\hdots,n\}$ around $0\in\mathbb R^d$ and we define a functional design $\{x_i, i=1,\hdots,n\}$ verifying
\begin{equation}\label{eq:gen_design}
x_i:=x_0+\sum_{j=1}^d x_{i,j}\varphi_j.
\end{equation}

The advantage of such a method is its flexibility: all multivariate designs and all basis of $\mathbb H$ can be used. Then, by choosing an appropriate design and an appropriate basis, we can generate designs satisfying some constraints defined by the context.

Remark that $\{\varphi_1,...,\varphi_d\}$ can be seen as predefined optimisation directions in the sense that the input $x$ will vary exclusively in the directions of the space $\text{span}\{\varphi_1,...,\varphi_d\}$. Therefore their choice  have a great influence on the precision of the results and has to be made carefully.

\subsection{Data-driven directions of optimisation}

If a training sample $\{(X_i,Y_i),i=1,\hdots,n\}$ is available, with $X_i\in\mathbb H$, for all $i$ and $\mathbb E[Y_i]=m(X_i)$,  it may be relevant to use the information of this sample to find a suitable basis. We consider here two method to generate data-driven basis:
\begin{itemize}
\item Principal Components which is the basis of $\mathbb H$ verifying 
\[\frac{1}{n}\sum_{i=1}^n\|X_i-\widehat\Pi_d X_i\|^2=\min_{\Pi_d }\left\{\frac{1}{n}\sum_{i=1}^n\|X_i-\Pi_d X_i\|^2\right\},\]
where $\widehat{\Pi}_d$ is the orthogonal projector on $\text{span}\{\varphi_1,\hdots,\varphi_d\}$, $\|\cdot\|$ is a norm on the space $\mathbb H$ and the minimum on the right-hand side is taken over all orthogonal projectors $\Pi_d$ on $d$-dimensional subspaces of $\mathbb H$. 
\item Partial Least Squares \citep{wold_soft_1975,preda_pls_2005} which permits to take into account the interaction between $X$ and $Y$. It is computed iteratively by the procedure described in \citet{delaigle_methodology_2012}.  For theoretical results on the PLS basis see \citet{delaigle_methodology_2012,blazere_PLS_2014}. For practical implementation, see Algorithm~\ref{algo:PLS}.
\end{itemize}
\begin{algorithm}\small
 \KwData{Training sample $\left\{(X_i,Y_i), i=1,\hdots,n\right\}$}
 Initialization : 
 \begin{equation*}
 X_i^{[0]}=X_i-\frac{1}{n}\sum_{j=1}^n X_j,\ Y_i^{[0]}=Y_i-\frac{1}{n}\sum_{j=1}^nY_j
 \end{equation*}
 \For{$j=1,\hdots,d$}{
   Estimate $\varphi_j$ by the empirical covariance of $X_i^{[j-1]}$ and $Y_i^{[j-1]}$: \[\varphi_j=\sum_{i=1}^nY_i^{[j-1]}X_i^{[j-1]}/\left\|\sum_{i=1}^nY_i^{[j-1]}X_i^{[j-1]}\right\|\]
   Fit the models $Y_i^{[j-1]}=\beta_j\langle X_i^{[j-1]},\varphi_j\rangle+\varepsilon_i^{[j]}$ and $X_i^{[j-1]}=\delta_j\langle X_i^{[j-1]},\varphi_j\rangle+W_i^{[j]}$ by least-squares that is 
   \[\widehat{\beta}_j:=\sum_{i=1}^nY_i^{[j-1]}\langle X_i^{[j-1]},\varphi_j\rangle/\sum_{i=1}^n\langle X_i^{[j-1]},\varphi_j\rangle^2\]
    and 
    \[\widehat{\delta}_j:=\sum_{i=1}^n\langle X_i^{[j-1]},\varphi_j\rangle X_i^{[j-1]}/\sum_{i=1}^n\langle X_i^{[j-1]},\varphi_j\rangle^2\]
   Define $X_i^{[j]}:=X_i^{[j-1]}-\langle X_i^{[j-1]},\varphi_j\rangle\widehat{\delta}_j$ and $Y_i^{[j]}:=Y_i^{[j-1]}-\widehat{\beta}_j\langle X_i^{[j-1]},\varphi_j\rangle$ the residuals of the two fitted models\;}
 \caption{\label{algo:PLS} practical implementation of PLS basis \citep[Section A.2]{delaigle_methodology_2012}}
\end{algorithm}

\subsection{Multivariate designs}
In this article, we focus on the most classical designs. However, the method we propose is flexible and can be used with any multivariate design, for instance Latin Hypercube Sampling \citep[see][and references therein]{liu_two_2015}, small composite designs~\citep{draper_small_1990}, augmented-pair designs \citep{morris_class_2000}... We also refer to \citet{georgiou_class_2014} and references therein for the recent advances on the subject. 

The \textit{$2^d$ factorial design} is one of the simplest. It is a \textit{first-order design} is the sense that it is frequently used to fit a first-order linear model. For each explanatory variable $ x_1, \hdots, x_d $, we choose two levels (coded by $+1$ and $-1$) and we take all the $ 2^d $ combinations of these two levels. When $ d $ is large, it may be impossible to achieve the $ 2^d $ factorial experiments, hence \textit{fractional factorial design} keeps only a certain proportion (e.g. a half, a quarter,...) of points of a $2^d$ factorial design. Typically, when a fraction $1/(2^p)$ is kept from the original $2^d$ design, this design is called $2^{d-p}$ factorial design.  The points removed are carefully chosen, we refer e.g. to \citet{gunst_fractional_2009} for more details. In our context, since we have the freedom to choose the dimension $d$, the interest of considering fractional factorial designs relies on its flexibility. For a given number $2^k$ of design points, all pairs $(d,p)$ of positive integers such that $d-p=k$  gives a different design, the choice $d=k$ and $p=0$ leads to the full factorial design while larger values of $p$ allow to explore a higher-dimensional space keeping the number of experiments low. 

Traditional second-order designs are factorial designs, central composite designs and Box-Behnken designs. 
\begin{itemize}
\item $3^d$ or $3^{d-p}$ factorial designs are similar to $2^d$ and $2^{d-p}$ factorial designs but with three levels ($+1$, $-1$ and $0$).
\item \textit{Central Composite Designs (CCD)} are obtained by adding to the two-level factorial design (fractional or not) two points on each axis of the control variables on both sides of the origin and at distance $\alpha>0 $ from the origin. 
\item \textit {Box-Behnken Designs (BBD)} are widely used in the industry. It is a well-chosen subset of the $3^d$ factorial design. Box-Behnken designs are not used when $d=2$. 
For $ d \geq  4$, we refer to \citet[7.4.7]{myers_response_2009}.
\end{itemize}

\subsection{Design properties}
\label{sec:lsquares_design_prop}

One of the interests of the design generation method\eqref{eq:gen_design} is that all design properties (orthogonality, rotatability and alphabetic optimality) verified by the multivariate design $\{\mathbf x_i, i=1,...,n\}$ are also verified for the design $\{x_i, i=1,...,n\}\subset\mathbb H$. 

To explain that we focus on first-order and second-order designs but the same reasoning may apply to other models and other kind of optimality properties related to the model considered. Let us first rewrite these models.
 
 \subsubsection*{First-order model}
Recall that, for all $i=1,\hdots,n$, $x_i=x_0+\sum_{j=1}^dx_{i,j}\varphi_j$, then the first-order model can be rewritten
  \begin{equation}\label{model2}
  y_i: = \alpha+\langle\beta, x_0\rangle+\sum_{j=1}^d x_{i,j}\langle\beta,\varphi_j\rangle + \varepsilon_i,\text{ for }i=1,\hdots,n.
  \end{equation}
  With our choice of design points, this model is a first-order multivariate model and can be written
  \begin{equation}\label{eq:linmodmatrix}
  \mathbf Y=\mathbf X\boldsymbol\beta+\boldsymbol\varepsilon
  \end{equation} 
  with design matrix
  \[\mathbf X=\begin{pmatrix}
  1 & x_{1,1} & \hdots & x_{1,d} \\ 
  1 & x_{2,1} &  & x_{2,d} \\ 
  \vdots & \vdots & \ddots & \vdots \\ 
  1 & x_{n,1} & \hdots & x_{n,d}
  \end{pmatrix} \]
  and coefficients $\boldsymbol\beta=(\alpha+\langle\beta, x_0\rangle,\langle\beta,\varphi_1\rangle,\hdots,\langle\beta,\varphi_d\rangle)^t$. Then, the first-order linear model in $\mathbb H$ is in fact a first-order multivariate linear model, with inputs $\{\mathbf x_i, i=1,...,n\}$.
\subsubsection*{Second-order model}
Now we can see that a similar conclusion holds for the second-order model, which can also be written $\mathbf Y=\mathbf X\boldsymbol\beta+\boldsymbol\varepsilon$ with 
\[\mathbf X=\begin{pmatrix}
  1 & x_{1,1} & \hdots & x_{1,d}& x_{1,1}^2 &x_{1,1}x_{1,2} & \hdots &  x_{1,d}^2\\ 
  1 & x_{2,1} &  & x_{2,d} & x_{2,1}^2 &x_{2,1}x_{2,2} & \hdots &  x_{2,d}^2\\ 
  \vdots & \vdots & \ddots & \vdots  &\vdots & \vdots & \ddots & \vdots \\ 
  1 & x_{n,1} & \hdots & x_{n,d} & x_{n,1}^2 &x_{n,1}x_{n,2} & \hdots &  x_{n,d}^2
  \end{pmatrix} \]

\subsubsection{Design properties}

An important property is the \textit{orthogonality}. 
An orthogonal design is a design for which the matrix $\mathbf X^t\mathbf X$ is diagonal. This implies that the vector $\widehat{\boldsymbol\beta}$ is also a Gaussian random vector with independent components and makes it easier to test the significance of the components of $\boldsymbol\beta$ in the model. $2^d$ factorial designs are orthogonal first-order designs. However, fractional designs have to be constructed carefully in order to keep the orthogonality property. For second-order designs, we refer to \citet{box_multi-factor_1957} for general criteria applied to factorial and fractional factorial designs. Central Composite Designs are orthogonal if $\alpha=\sqrt{\left(\sqrt{F(F+2d+n_0)}-F\right)/2}$,
where $F$ is the number of points of the initial factorial design (see \citealt{myers_response_2009}).

A design is said to be \textit{rotatable} if $\text{Var}(\widehat{y}(\mathbf x))$ depends only on the distance between $\mathbf x$ and the origin. This implies that the prediction variance is unchanged under any rotation of the coordinate axes. We refer to \citet{box_multi-factor_1957} for conditions of rotatability. All first-order orthogonal designs are also rotatable. This is not the case for second-order designs, for instance a CCD design is rotatable if $\alpha=F^{1/4}$ which means that a CCD design can be rotatable and orthogonal only for some specific values of $n_0$ and $F$. Box-Behnken designs are rotatable for $d=4$ and $d=7$. Some measures of rotatability have been introduced \citep{khuri_measure_1988,draper_index_1988,draper_another_1990,park_measure_1993} in order to measure how close a design is to the rotatability property. 


  The important point is that all the design properties cited above only depends on the design matrix $\mathbf X$. Hence, all properties of the multivariate design $\{\mathbf x_i, i=1,...,n\}$ are automatically verified for the design $\{x_i, i=1,...,n\}$. 

Since all alphabetic optimality criteria~\citep{pazman-foundations_1986} are also exclusively based on properties of the design matrix $\mathbf X$, it is possible to define e.g. $D$-optimal designs for data in $\mathbb H$ with Equation~\eqref{eq:gen_design} by taking a $D$-optimal multivariate design.

  \section{Numerical experiments}
  
 \label{sec:num_exp}
In this section, $\mathbb H=\Ld$. 

\subsection{Functional designs}
\label{sec:func_design}

We use here the functions \textit{cube}, \textit{ccd} and \textit{bbd} of the R-package \textit{rsm} \citep{rsm} to generate respectively $2^d$ factorial designs, Central Composite Designs (CCD) and Box-Behnken Designs (BBD). 

\subsubsection*{Functional designs with Fourier basis}
In this section, we set $\varphi_1\equiv 1$ and for all $j\geq 1$, for all $t\in[0,1]$,
\[\varphi_{2j}(t)=\sqrt{2}\cos(2\pi jt)\text{ and } \varphi_{2j+1}(t)=\sqrt{2}\sin(2\pi jt).\]

\begin{figure}
\begin{tabular}{ccc}
Factorial $2^d$ design & CCD & BBD\\
\includegraphics[width=0.3\textwidth]{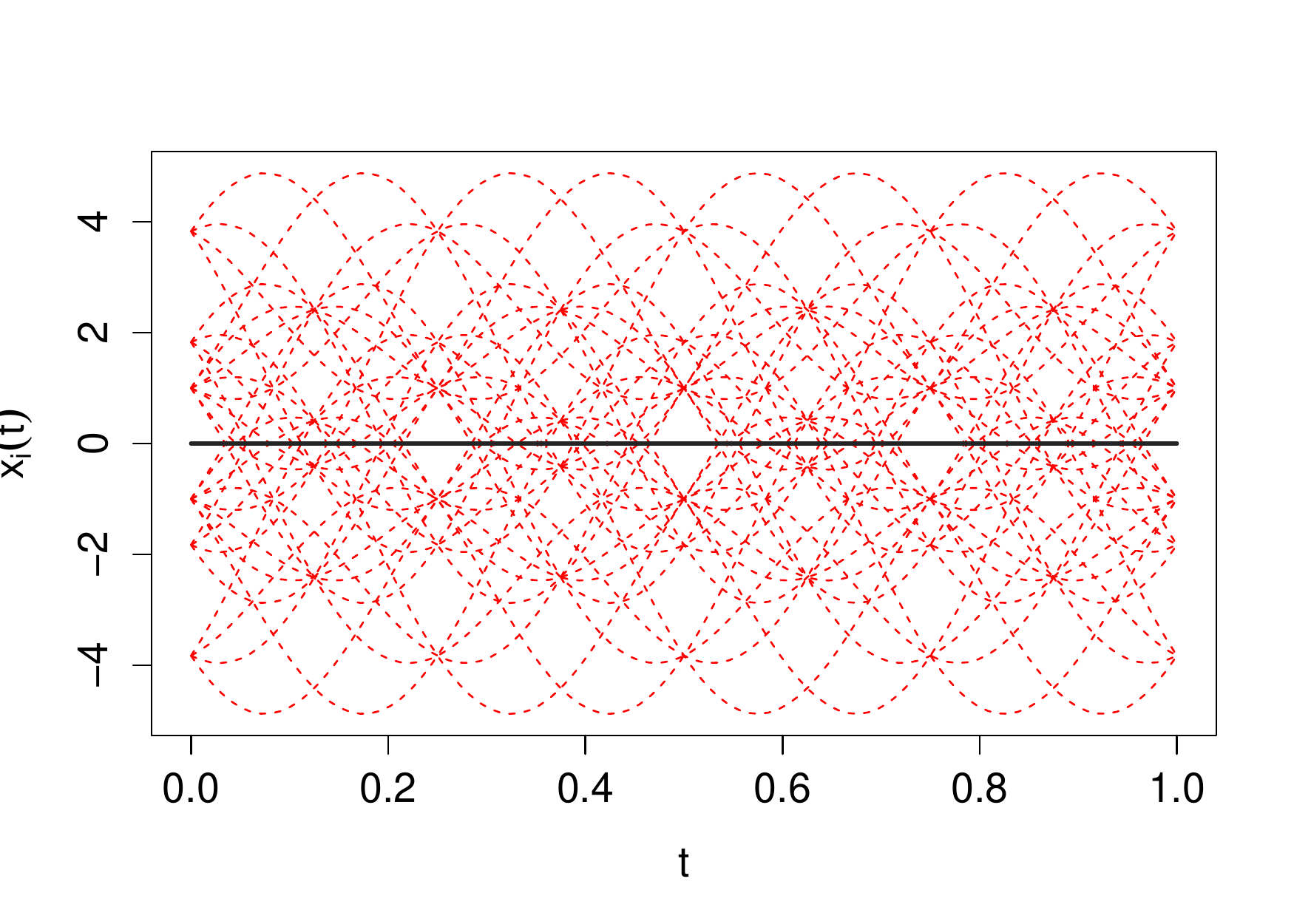}&\includegraphics[width=0.3\textwidth]{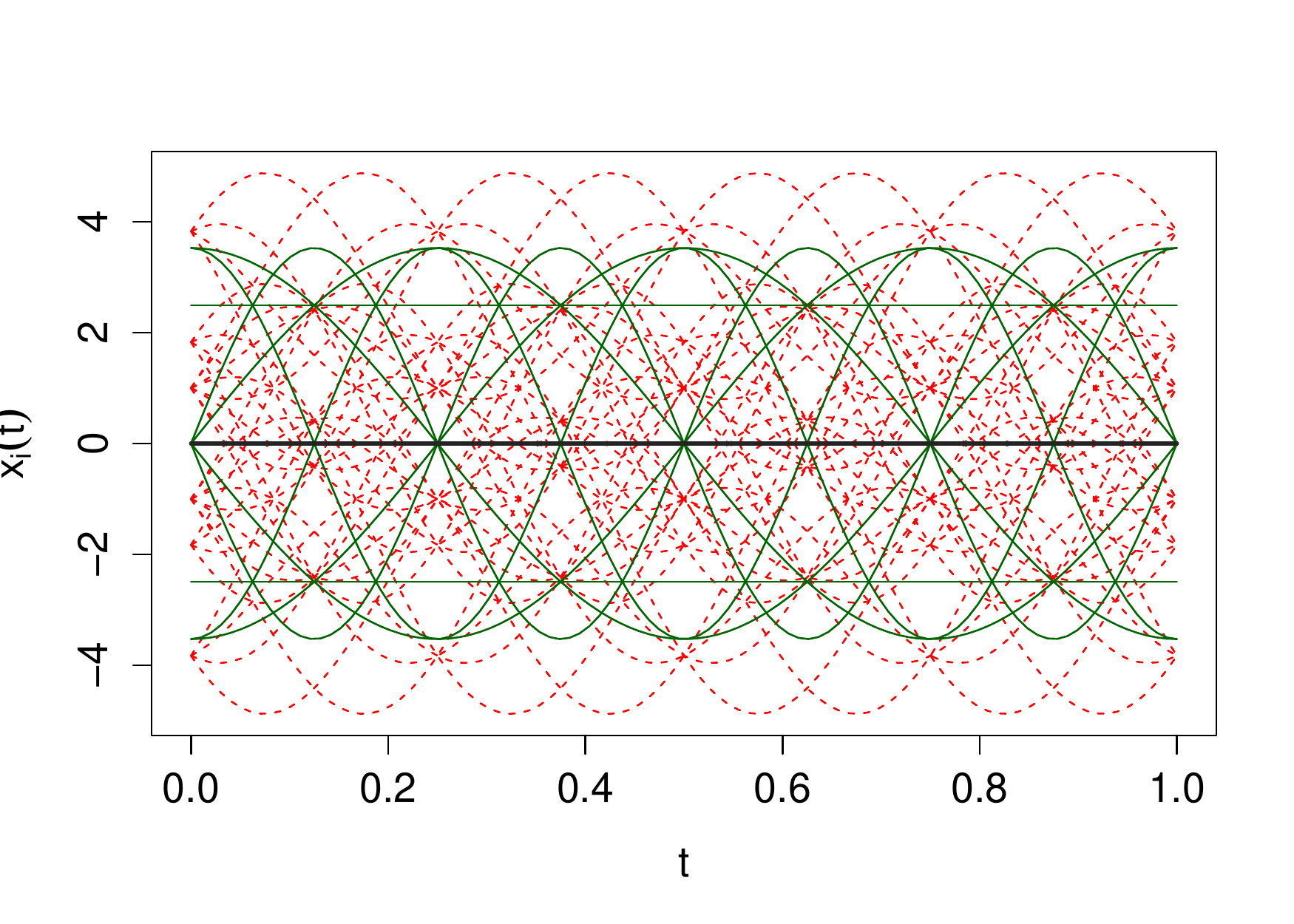}&\includegraphics[width=0.3\textwidth]{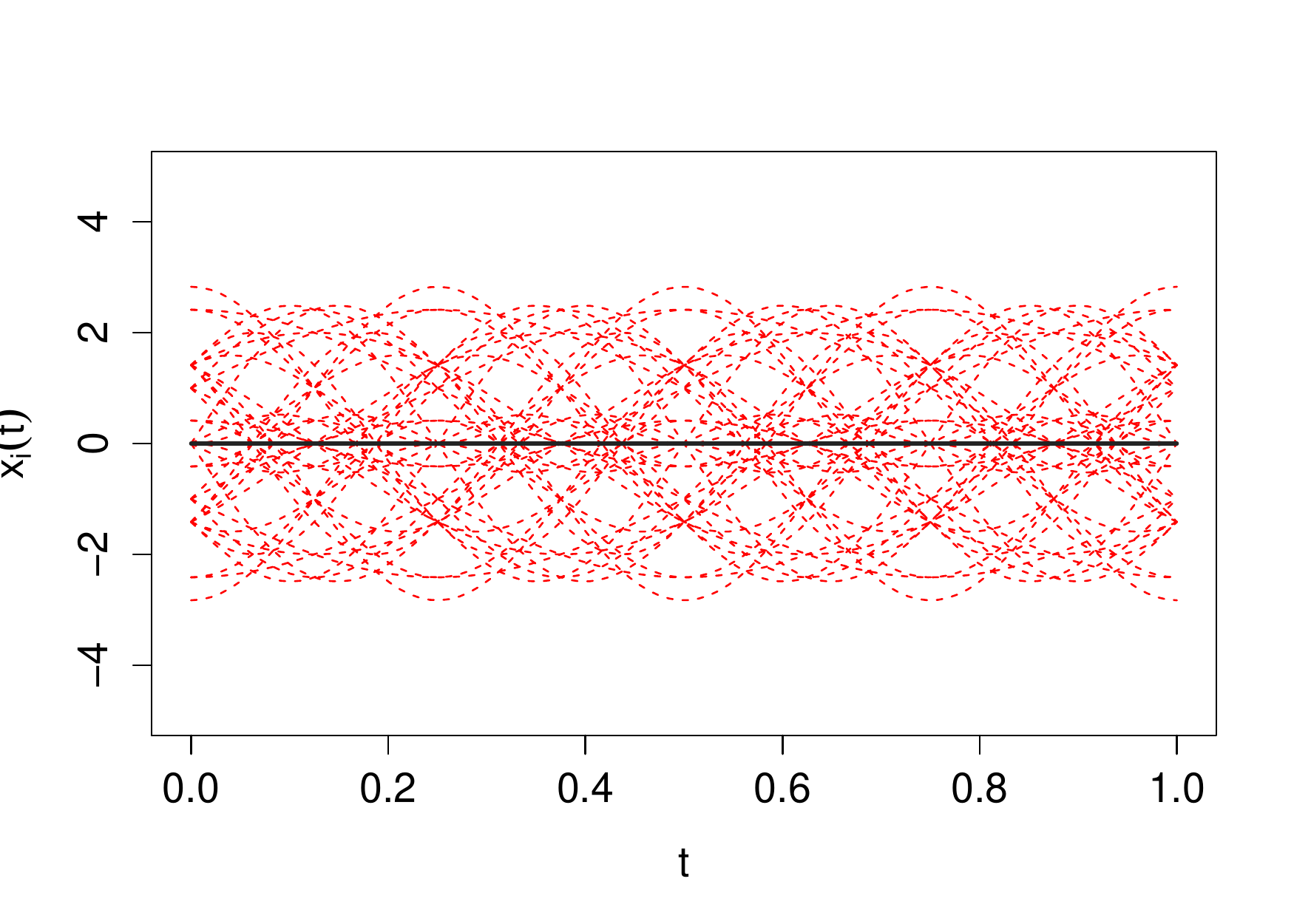}
\end{tabular}
\caption{\label{fig:fonctional_fourier}Functional designs with the Fourier basis ($d=5$). Gray thick line: $x_0\equiv 0$, red lines: points of the original $2^d$ or $3^d$ (for BBD) factorial design, green dotted lines: points added to the factorial design (for CCD). }
\end{figure}

The curves of the generated designs are given in Figure~\ref{fig:fonctional_fourier}. 

\subsubsection*{Functional design with data-driven bases}

We simulate a sample $\{X_1,\hdots,X_n\}$ comprising of $n=500$ realizations of the random variable
\[X(t)=\sum_{j=1}^J\sqrt{\lambda_j}\xi_j\psi_j(t),\]
with $J=50$, $\lambda_j=e^{-j}$, $(\xi_j)_{j=1,\hdots,J}$ an i.i.d. sequence of standard normal random variables and $\psi_j(t):=\sqrt{2}\sin(\pi(j-0.5)t)$. 

\begin{figure}
\begin{tabular}{ccc}
Factorial $2^d$ design & CCD & BBD\\
\includegraphics[width=0.3\textwidth]{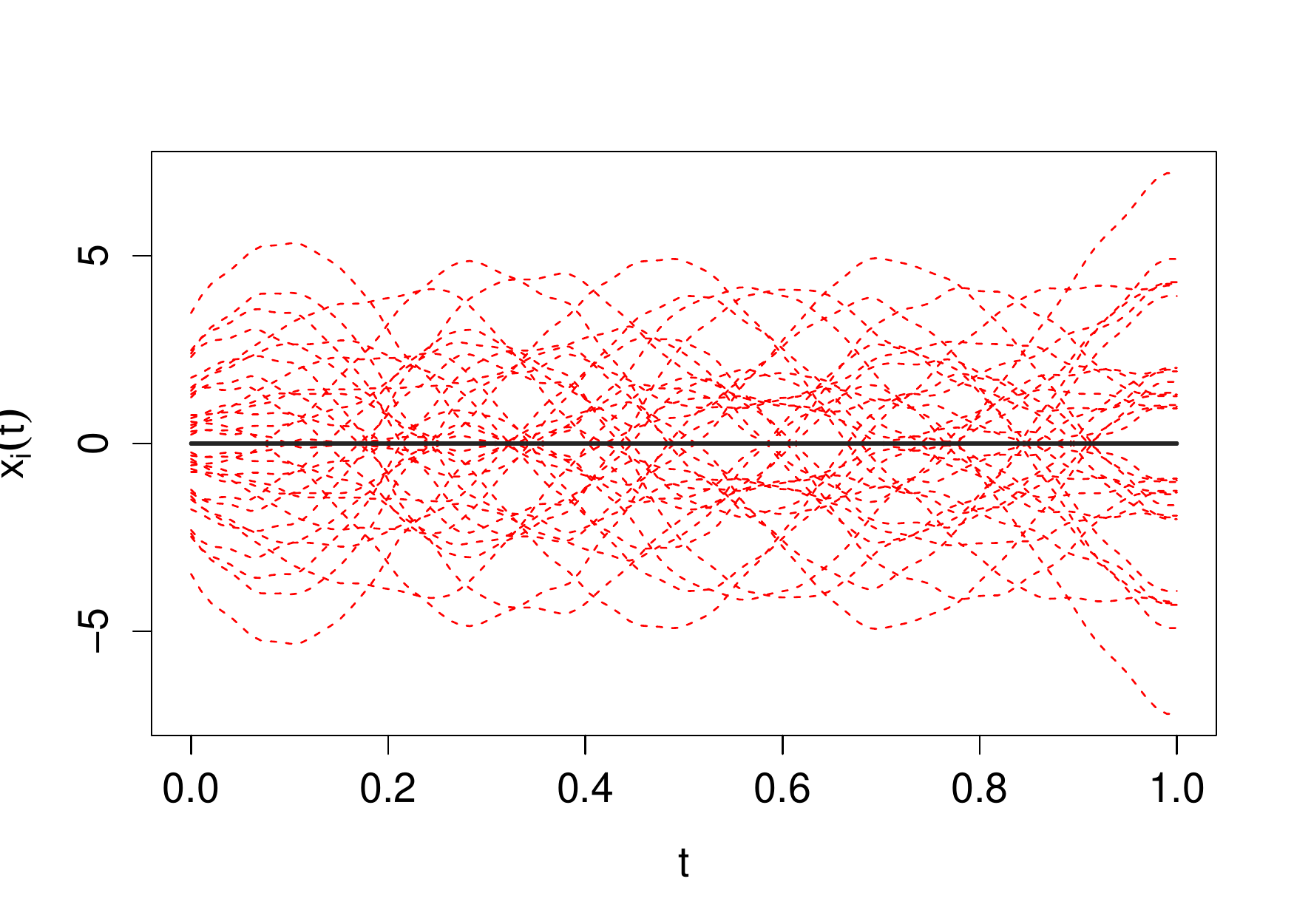}&\includegraphics[width=0.3\textwidth]{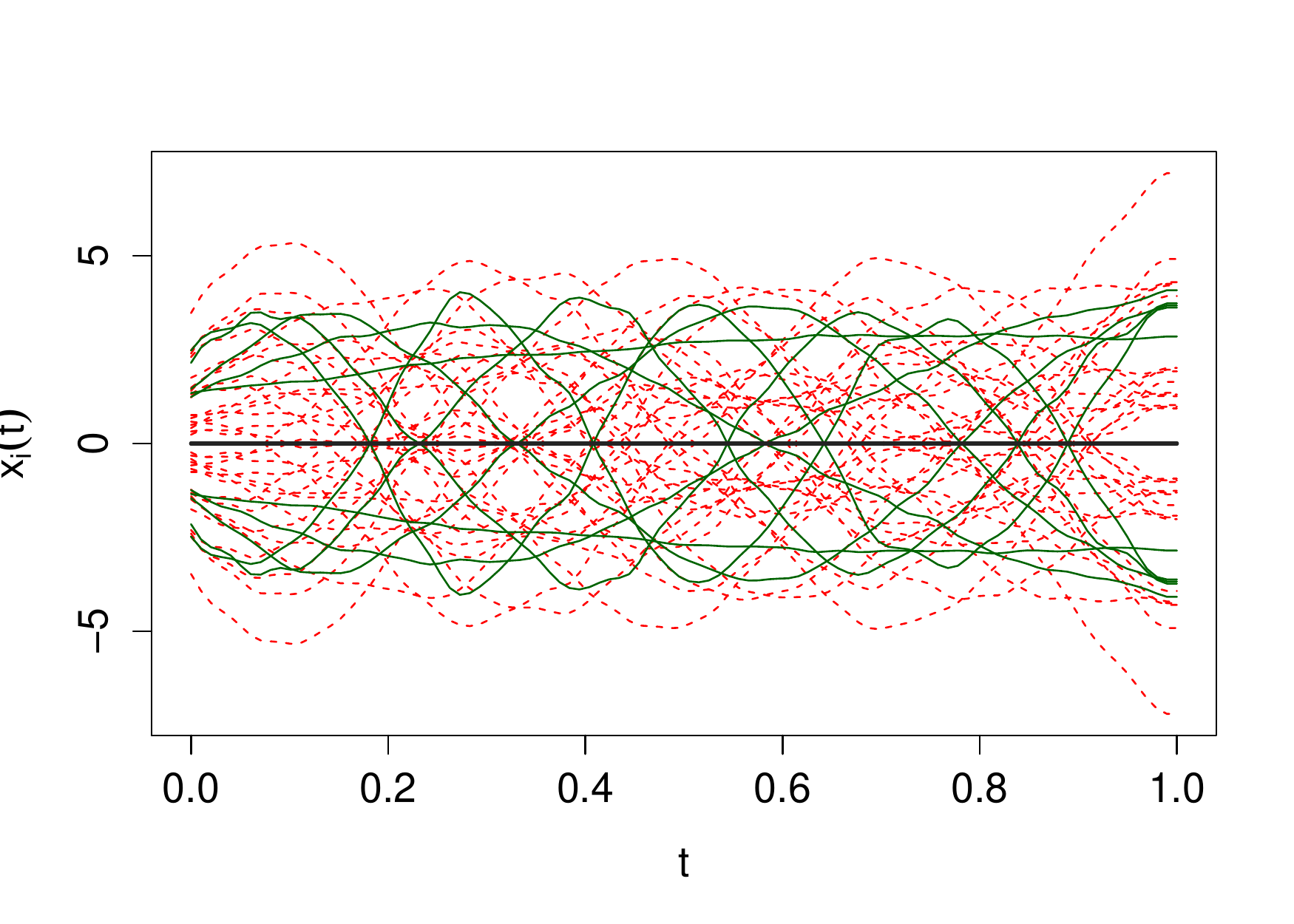}&\includegraphics[width=0.3\textwidth]{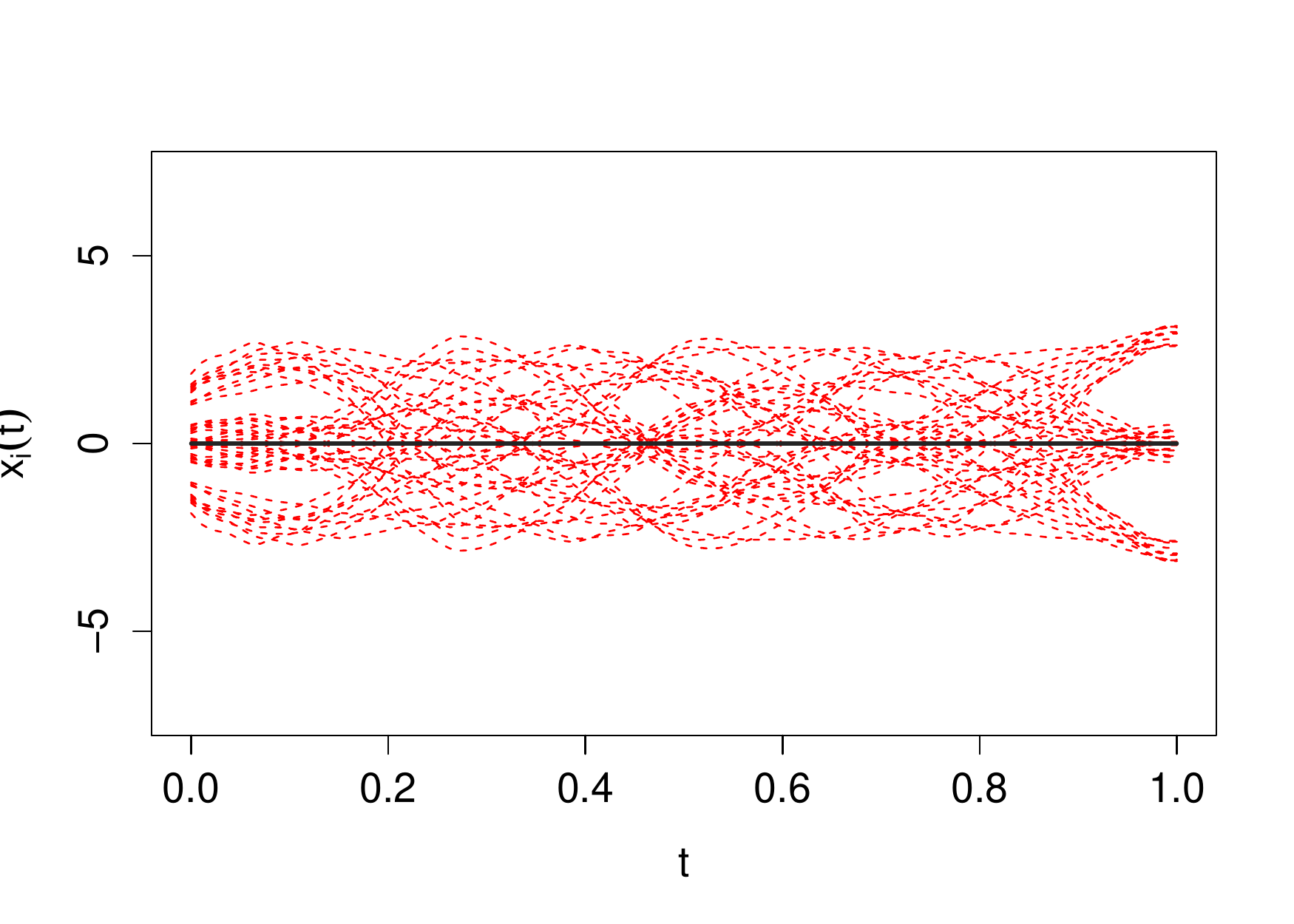}
\end{tabular}
\caption{\label{fig:fonctional_ACP}Functional designs with the PCA basis associated to $\{X_i, i=1,\hdots,n\}$ ($d=5$). The legend is the same as the one of Figure~\ref{fig:fonctional_fourier}.}
\end{figure}

The PCA basis only depends on  $\{X_i, i=1,\hdots,n\}$. In order to see the influence of the law of $Y$ on the PLS basis we define two training samples $ \left\{(X_i,Y_i^{(j)}), i=1,\hdots,n\right\}$ for $j=1,2$ with 
\[Y_i^{(j)}:=m_j(X_i)+\varepsilon_i ,\]
 $m_j(x):=\|x-f_j\|^2$, where
 \begin{eqnarray*}
 f_1(t)&:=&\cos(4\pi t)+3\sin(\pi t)+10,\\
 f_2(t)&:=&\cos(8.5\pi t)\ln(4t^2+10)
 \end{eqnarray*}
and $\varepsilon_1,\hdots,\varepsilon_n$, i.i.d. $\sim\mathcal{N}(0,0.01)$. 

\begin{figure}
\begin{tabular}{ccc}
Factorial $2^d$ design & CCD & BBD\\
\includegraphics[width=0.3\textwidth]{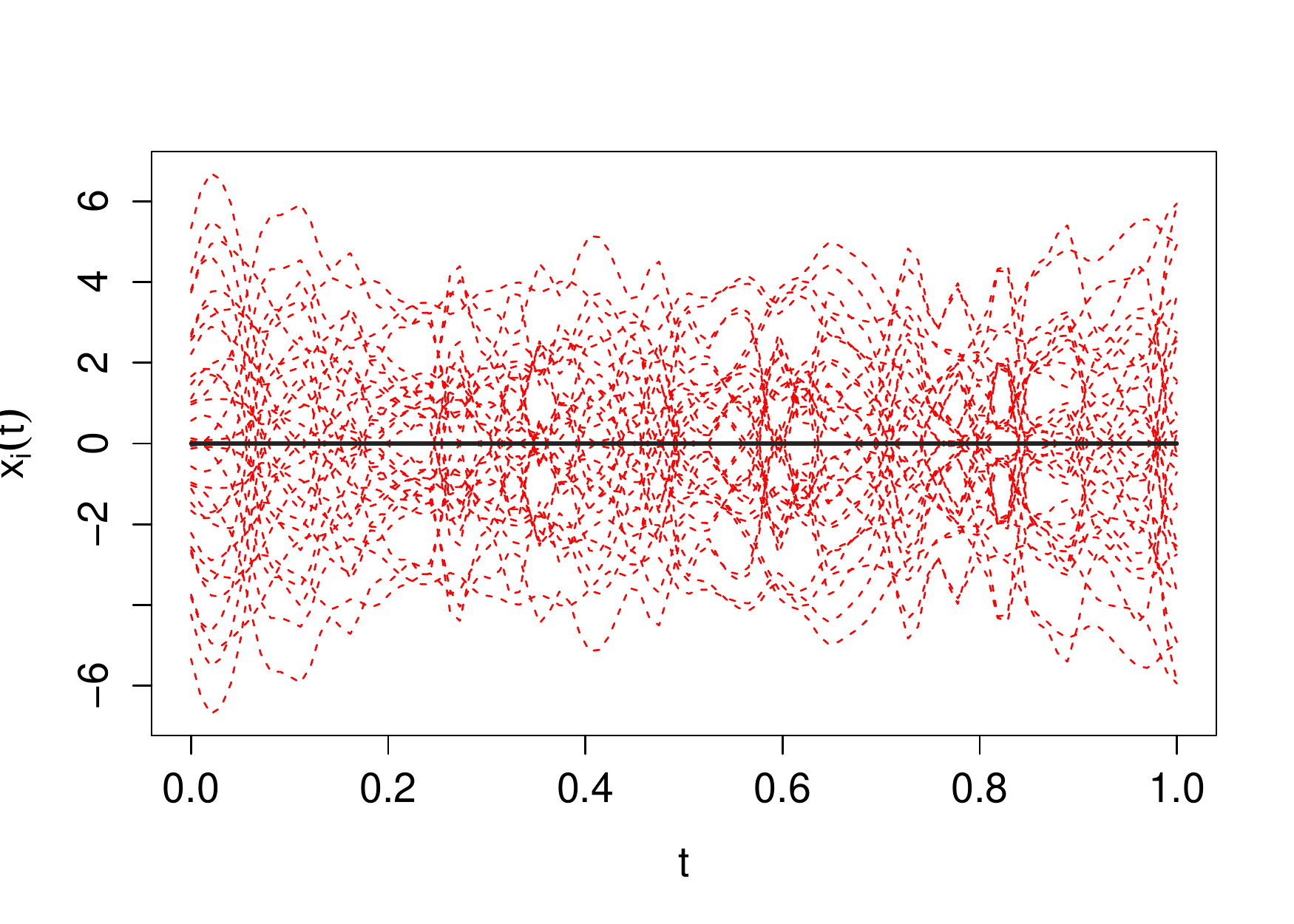}&\includegraphics[width=0.3\textwidth]{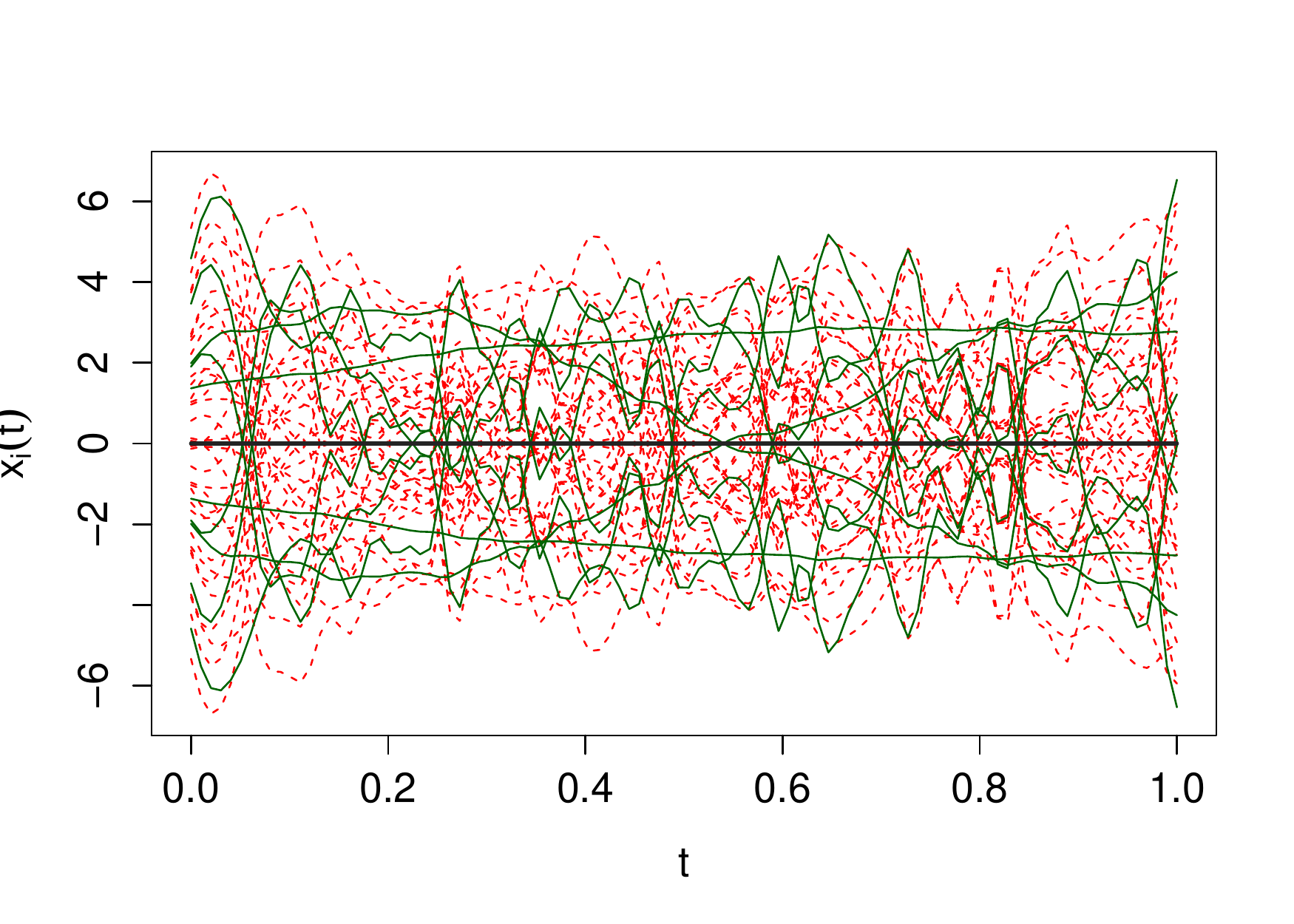}&\includegraphics[width=0.3\textwidth]{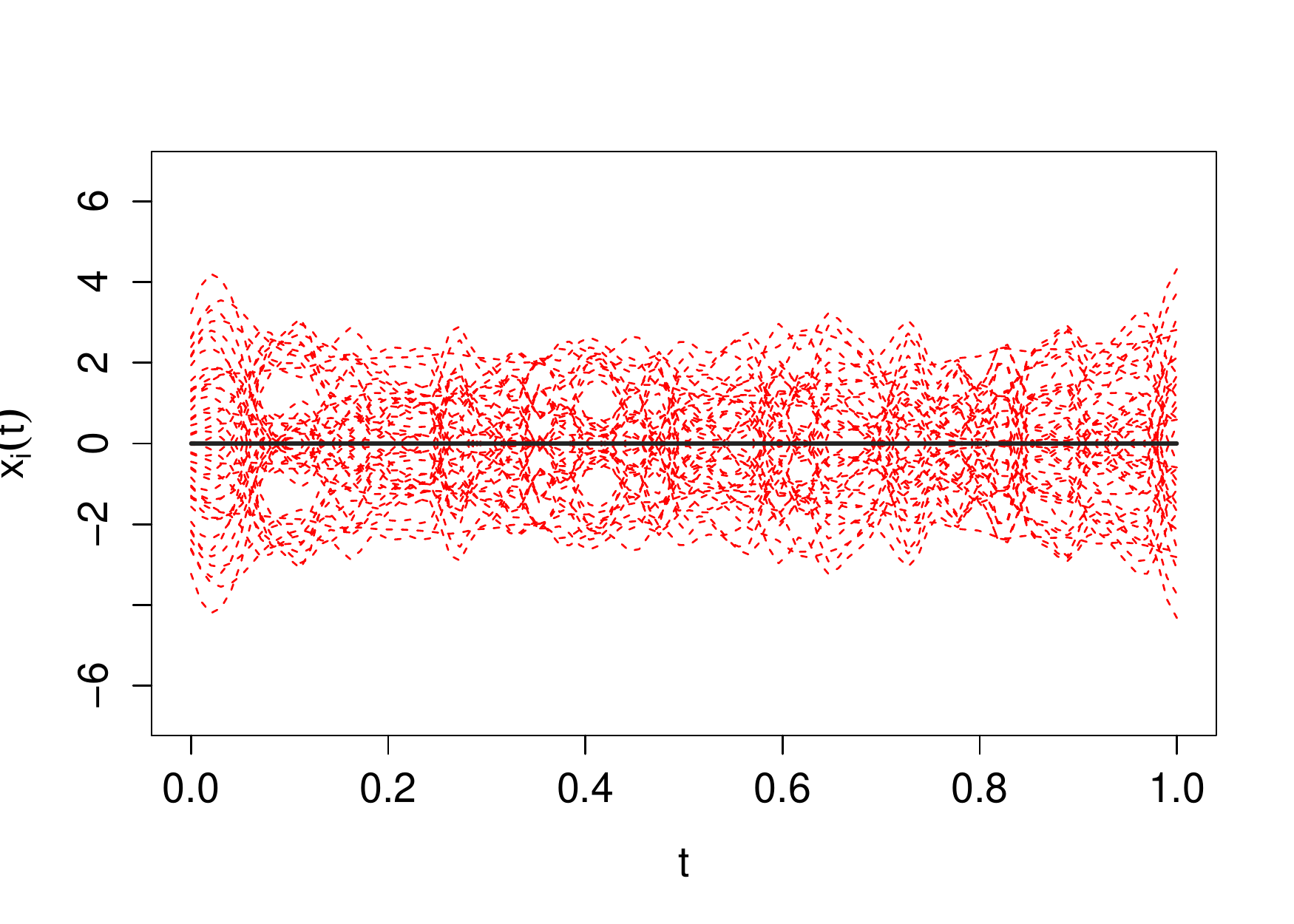}\\
\includegraphics[width=0.3\textwidth]{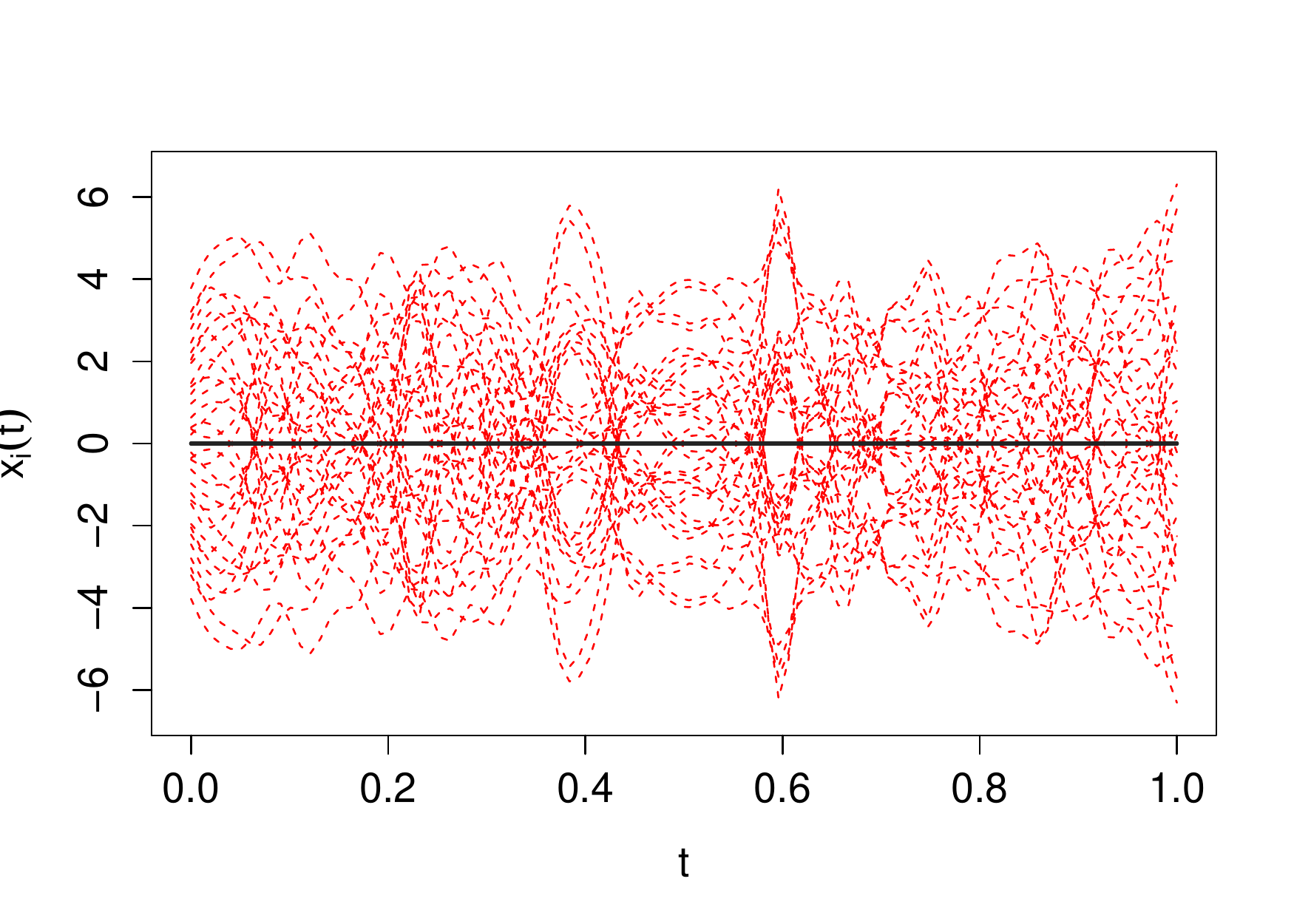}&\includegraphics[width=0.3\textwidth]{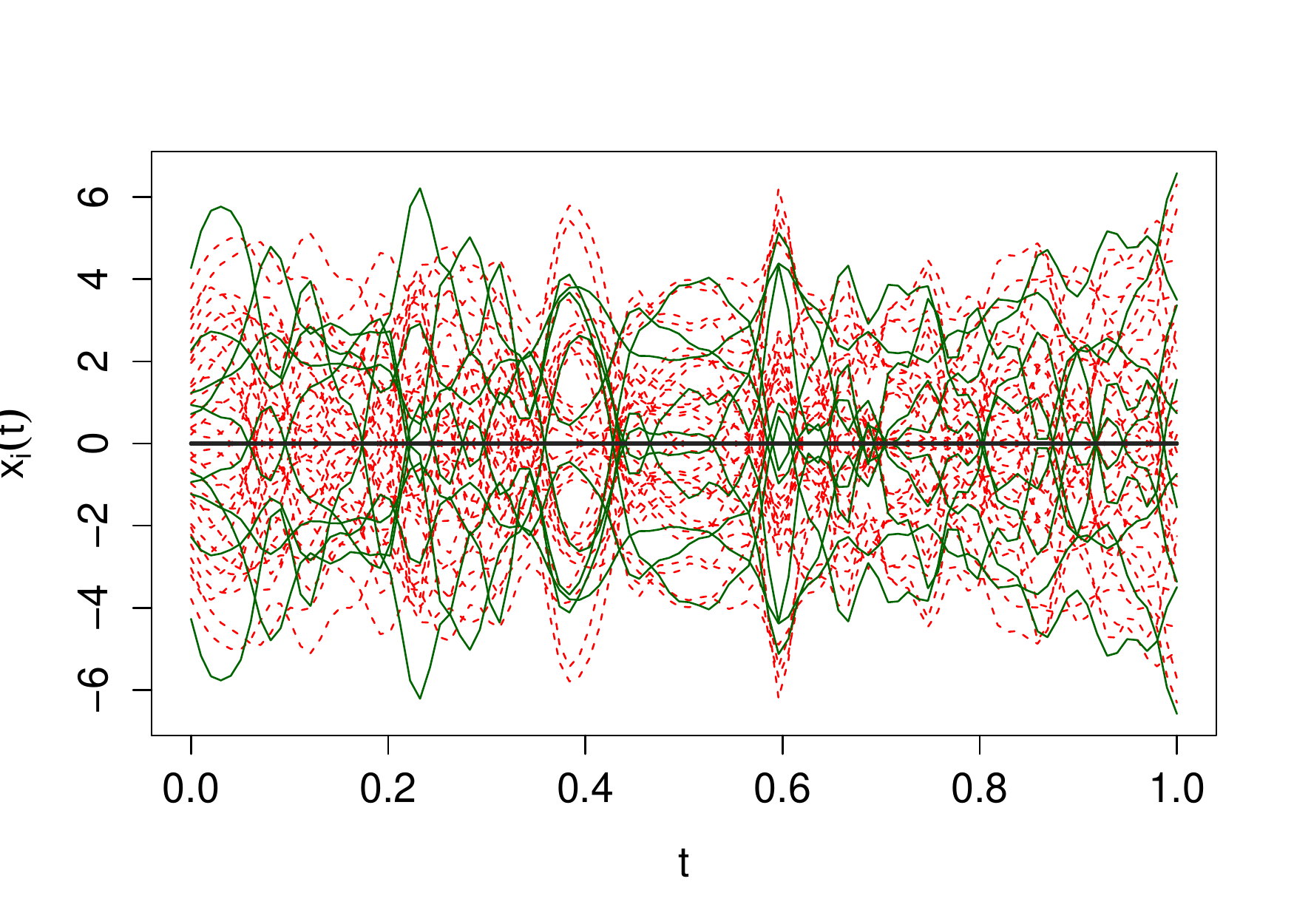}&\includegraphics[width=0.3\textwidth]{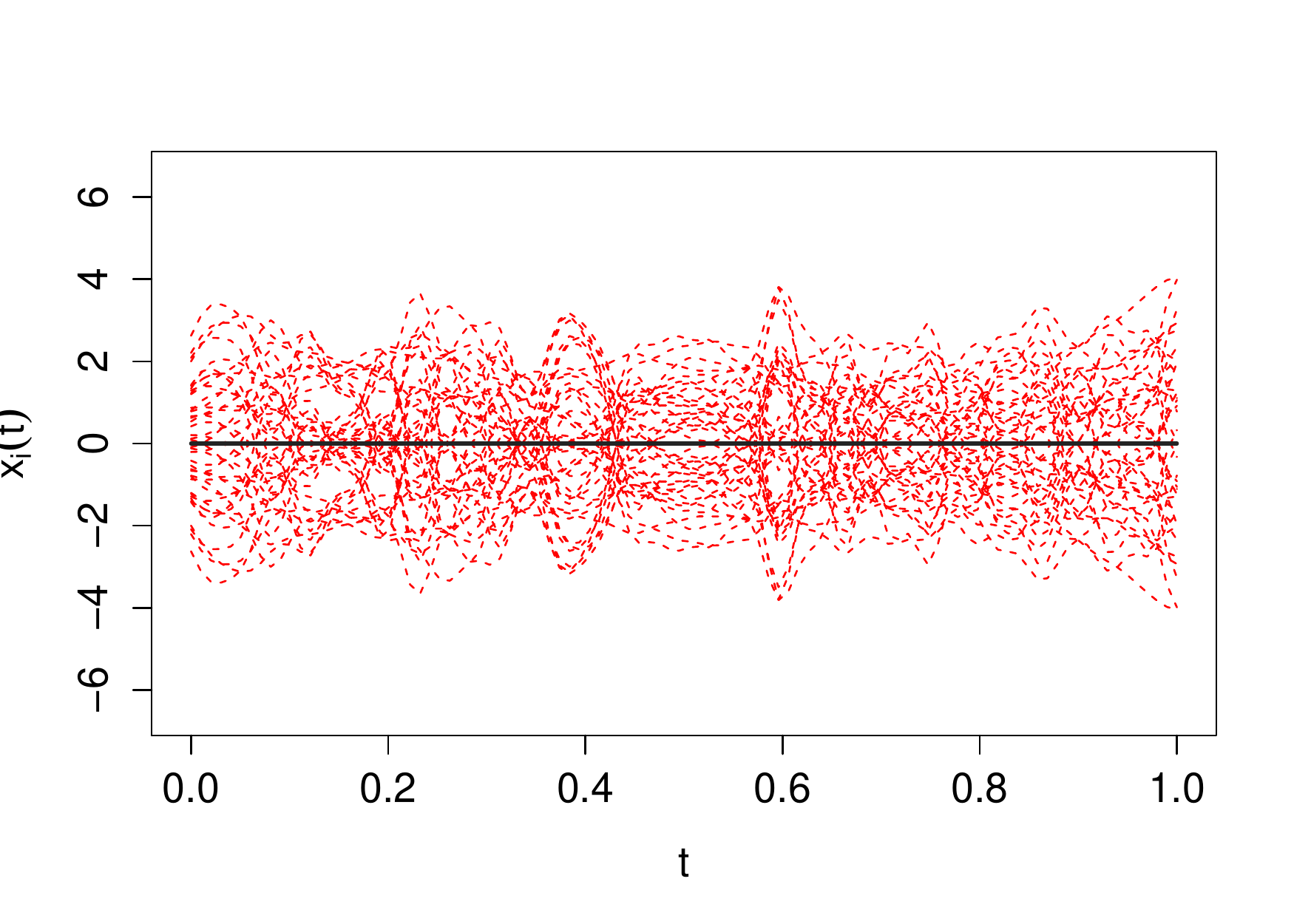}
\end{tabular}
\caption{\label{fig:fonctional_PLS}Functional designs with the PLS basis of the training sample $\{(X_i,Y_i^{(j)}, i=1,\hdots,n\}$, $j=1$ (first line) and $j=2$ (second line), $d=5$. The legend is the same as the one of Figure~\ref{fig:fonctional_fourier}.}
\end{figure}

The curves of the design generated by the PLS basis (Figure~\ref{fig:fonctional_PLS}) are much more irregular than those generated by the PCA basis (Figure~\ref{fig:fonctional_ACP}). However, remark that the designs generated by the PLS basis (Figure~\ref{fig:fonctional_PLS}) of the two samples show significant differences, which illustrates that the PLS basis effectively adapts to the law of $Y$. 
\subsection{Estimation of the response surface}
We use here the PLS basis calculated from the training sample $\{(X_i,Y_i^{(j)}),i=1,\hdots,n\}$ with $j=1$ or $j=2$. The aim is to approach the minimum $f_j$ of $m_j:x\in\mathbb H\to \|x-f_j\|^2$. 

We use the training sample a second time to determine the starting point of the algorithm. We take
\[x_0^{(0)}:=X_{i_{\min}}\text{ where }i_{\min}:={\arg\min}_{i=1,\hdots,n}\{Y_i\}.\]
The dimension is set to $d=8$. 

\subsubsection*{Approximation of $f_1=\cos(4\pi t)+3\sin(\pi t)+10$.}

\begin{description}
\item[Descent step] We generate a factorial $2^d$ design (Figure~\ref{fig:fonctional_PLS} -- left) $(x_1^{(0)},\hdots,x_{n_0}^{(0)})$ (here $n_0=2^d$) and we fit a first-order model
\[Y_i^{(0)}= \alpha^{(0)}+\sum_{j=1}^d\beta_j^{(0)}x_{i,j}^{(0)} +\varepsilon_i^{(0)},\]
to estimate the gradient. We realize two series of experiments along the line of steepest descent $x_0^{(0)}-\lambda_0\widehat\beta^{(0)}$ ($\lambda_0>0$). The first one (Figure~\ref{fig:dire_steep1}-- top left) suggests that the optimal value of $\lambda_0$ is between 0.4 and 0.6 and with the results of the second one we fix  $\lambda_0=0.50$. We set $x_0^{(1)}:=x_0^{(0)}-\lambda_0\widehat\beta^{(0)}$.

\begin{figure}
\begin{tabular}{cc}
\includegraphics[width=0.45\textwidth]{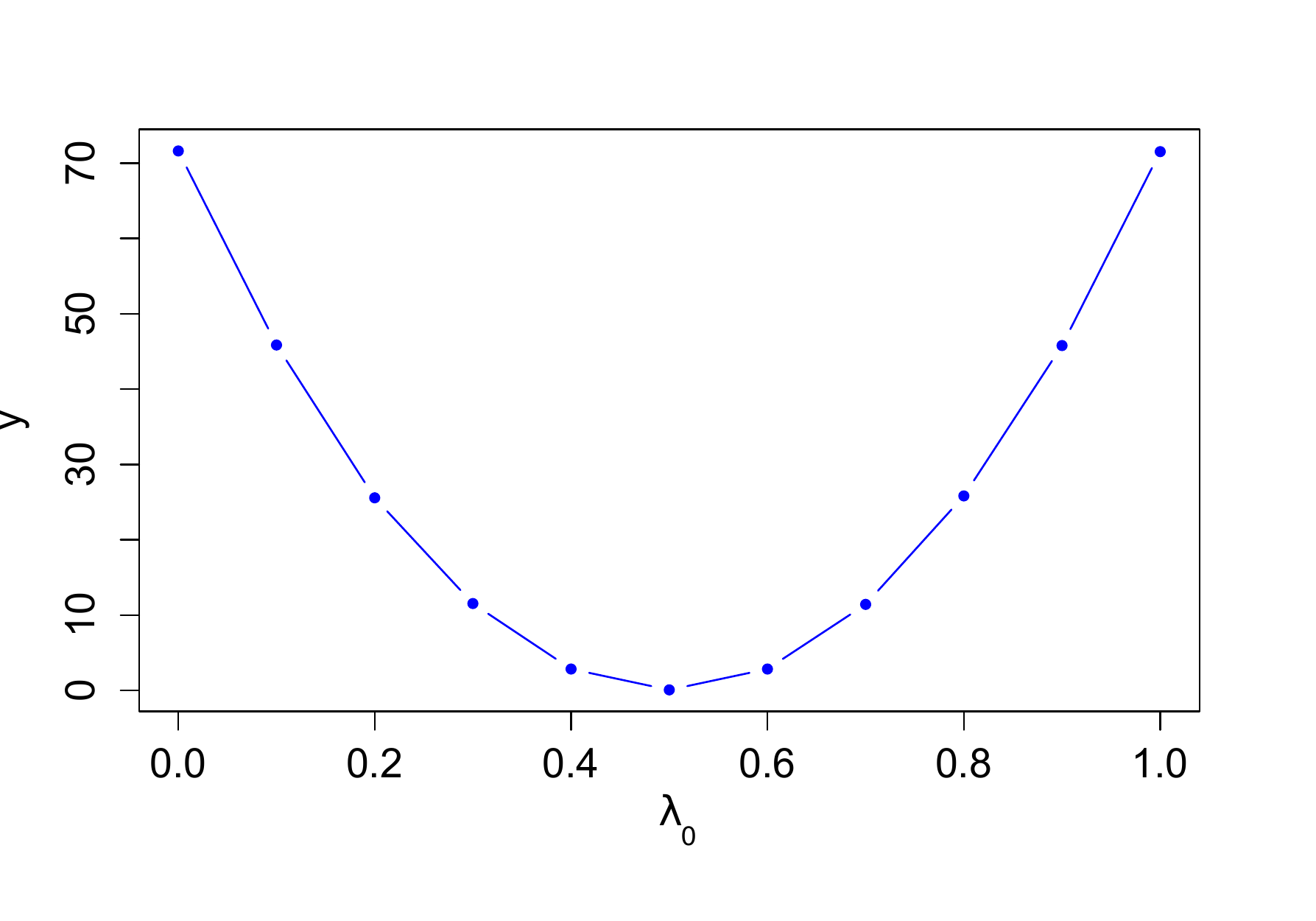}&\includegraphics[width=0.45\textwidth]{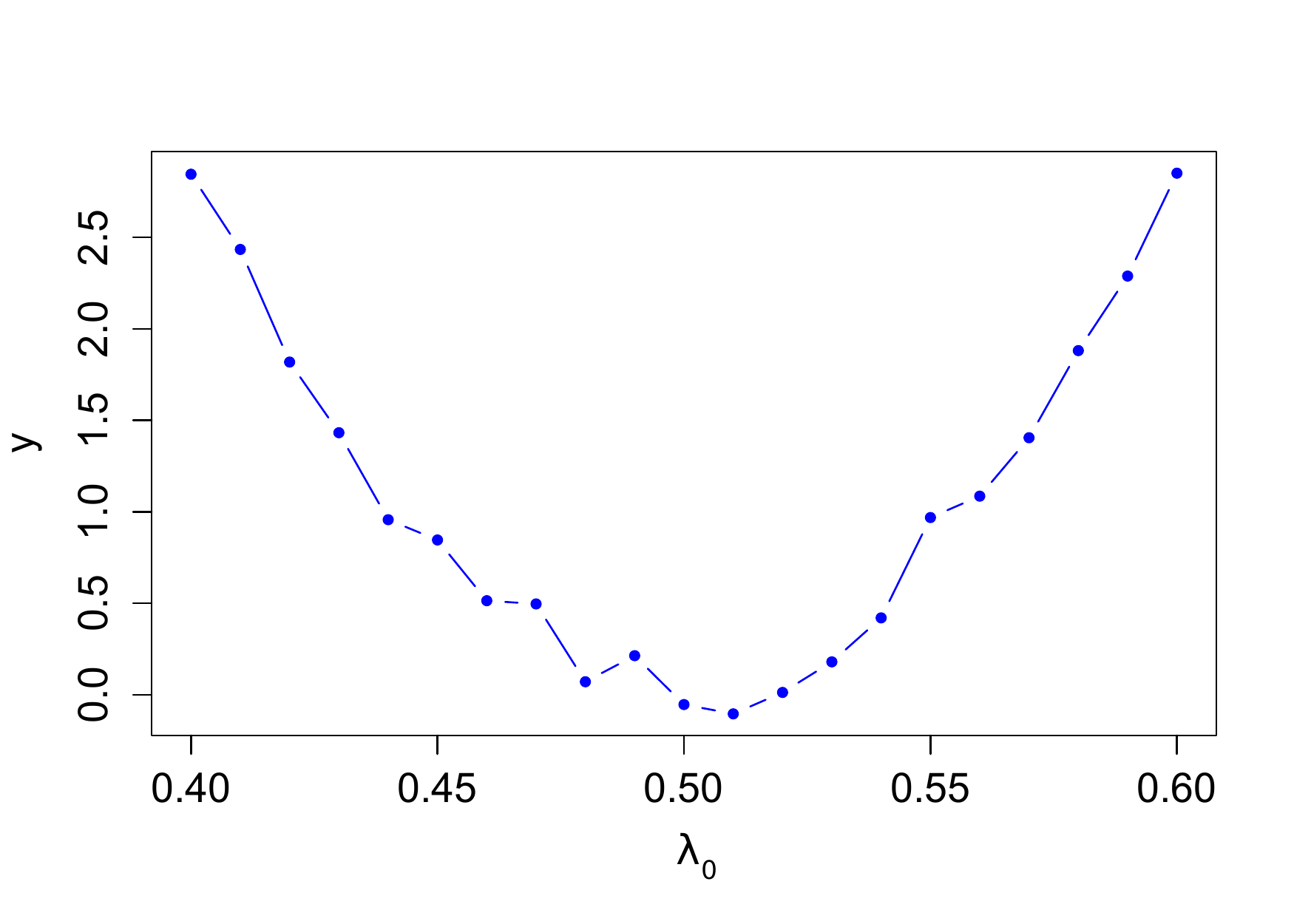}
\end{tabular}
\caption{\label{fig:dire_steep1} Results of experiments on the direction of steepest descent for the estimation of $f_1$. $x$-axis: $\lambda_0$, $y$-axis: response $Y=m_1(x_0^{(0)}-\lambda_0\widehat\beta^{(0)})+\varepsilon$. }
\end{figure} 

The value of $m_1$ at the starting point was $m(x_0^{(0)})=70.4\pm0.1$. At this step, we have $m_1(x_0^{(1)})=6.33\times 10^{-3}\pm 10^{-5}$ and we have done only $2^d+24=280$ experiments to reach this result.

We fit a first-order model once again with a $2^d$ factorial design and find that the norm of $\widehat\beta^{(1)}$ is very small ($\|\widehat\beta^{(1)}\|<0.02$) compared to $\|\widehat\beta^{(0)}\|=16.8\pm 0.1$ which suggests that we are very close to a stationary point. We also note that the $p$-value of the Fisher's test $H_0:\beta_1^{(1)}=\hdots=\beta_d^{(1)}=0$ against $H_1:\exists j\in\{1,\hdots,d\}, \beta_j^{(1)}\neq 0$ is very close to 1 which tends to confirm this assertion. 

\item[Final step] To improve the approximation, we fit a second-order model on the design points given by a Central Composite Design (Figure~\ref{fig:fonctional_PLS} -- center). The matrix $\widehat{H}$ at this step is an estimation of the matrix of the restriction to the space $\text{span}\{\varphi_1,\hdots,\varphi_d\}$ of the Hessian operator of $m$ at the point $x_0^{(1)}$. All the eigenvalues of $\widehat{H}$ are greater than $1.96>0$, this suggests that we are close to a minimum. We set $x_0^{(2)}:=-\widehat{H}^{-1}\widehat\beta^{(1)}$ and we have $m_1(x_0^{(2)}):=5.45\times 10^{-3}\pm 10^{-5}$. The CCD with $d=8$ counts 280 elements then we have realized $280$ experiments for the descent step plus $280$ for the final step, this rises to $557$ the total number of experiments performed. Figure~\ref{fig:result1} represents the different results. 

\begin{figure}
\centering
\includegraphics[width=0.6\textwidth]{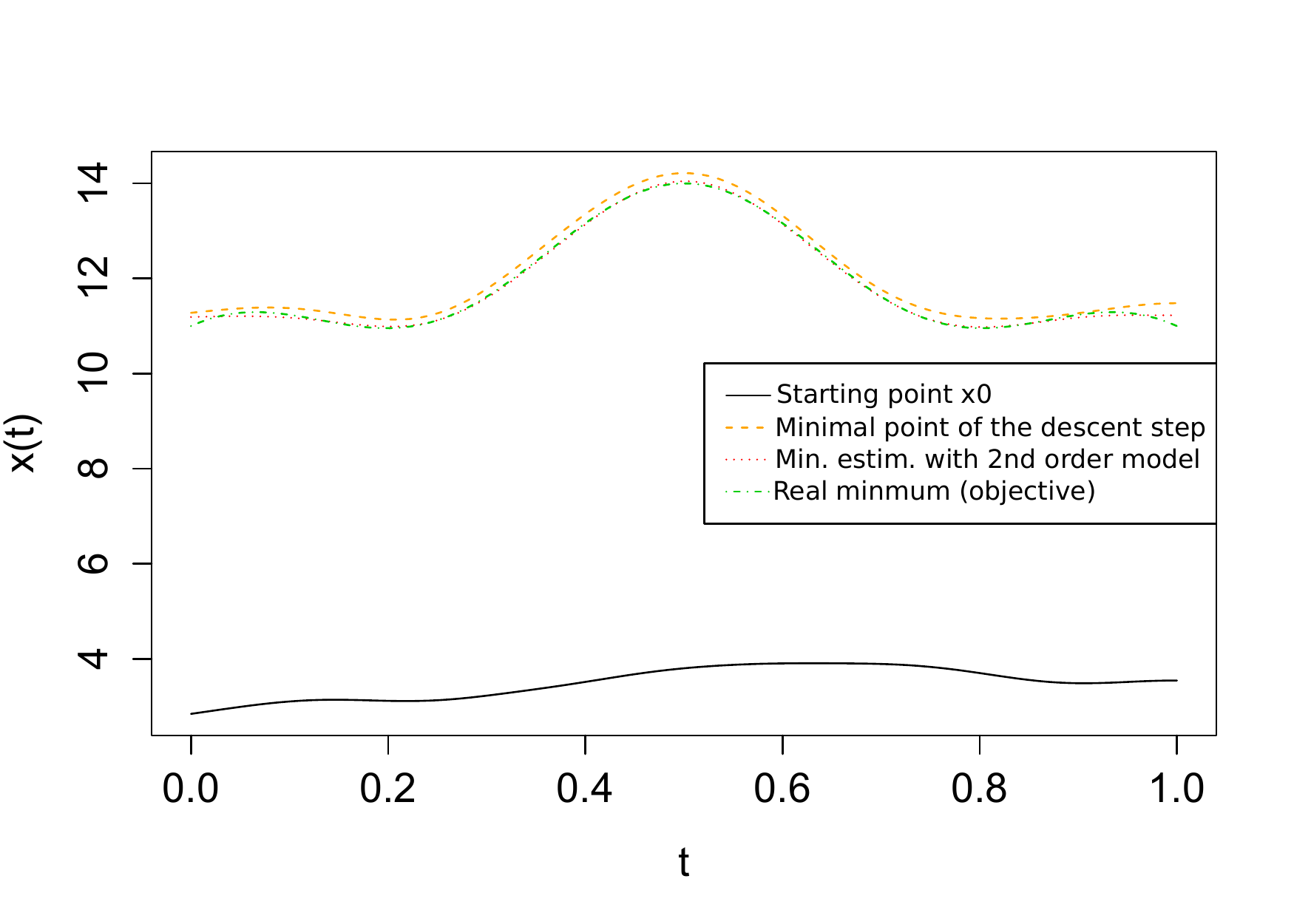}
\caption{\label{fig:result1} Result of optimization algorithm. }
\end{figure}
\end{description}

\subsubsection*{Approximation of $f_2(t)=\cos(8.5\pi t)\ln(4t^2+10)$}
We have here  $m_2(x_0^{(0)})=2.88\pm 0.01$. 

\begin{figure}
\centering
\includegraphics[width=0.5\textwidth]{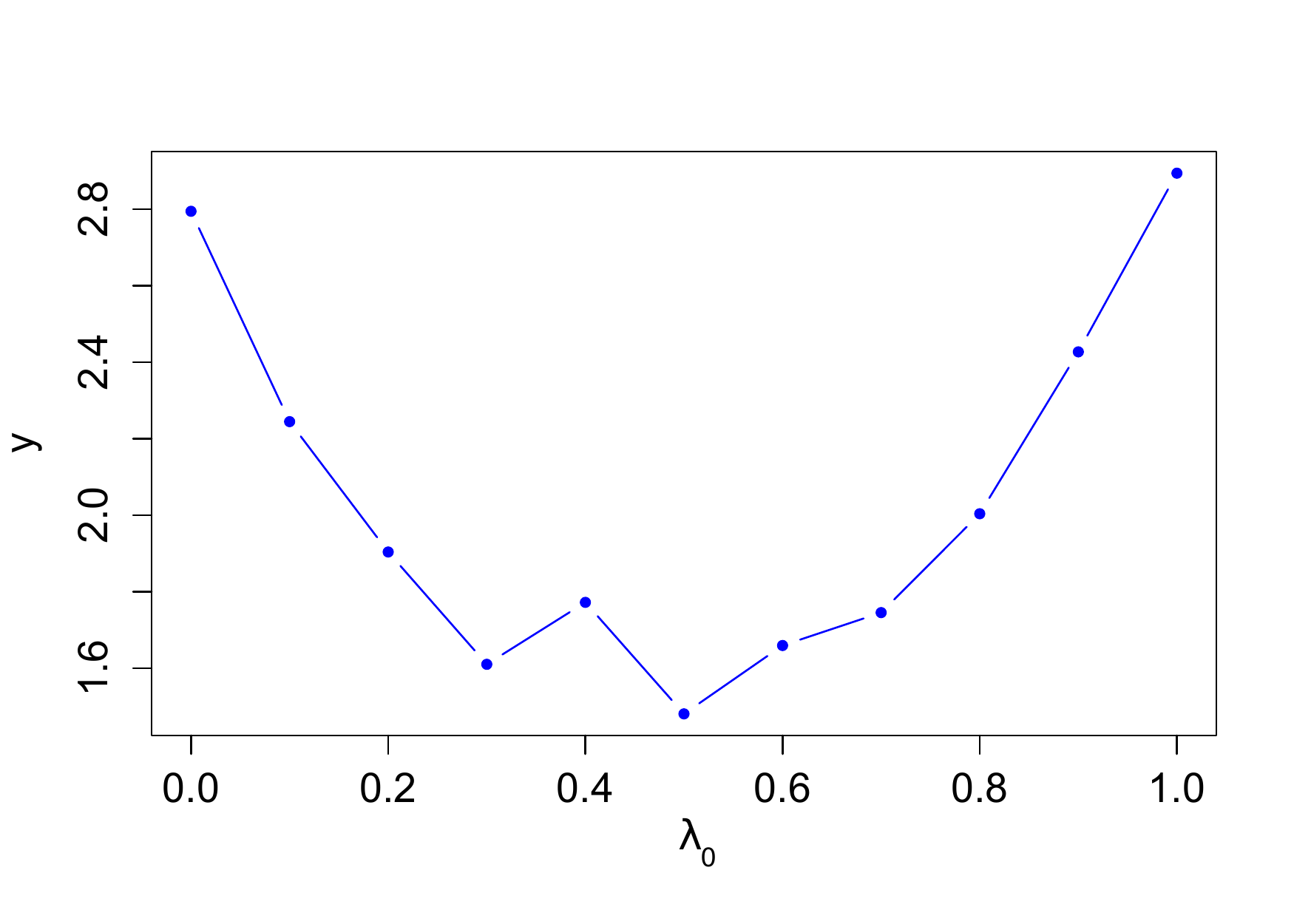}
\caption{\label{fig:dire_step2}Results of experiments on the direction of steepest descent for the estimation of $f_2$. 
}
\end{figure}

We follow the same steps as in the previous paragraph. Figure~\ref{fig:dire_step2}--left represents the evolution of the response along the direction of steepest descent. Here, since the response is noisy, refining the result without doing a too large number of experiments  seems to be difficult. Then, we fix $\lambda_0=0.5$ and $x_0^{(1)}=x_0^{(0)}-\lambda_0\widehat\beta^{(0)}$. We have $m_2(x_0^{(1)})=1.99\pm 0.01$. At this step, we have improved the response of about 31\%. This is not as important as the improvement of the first step of estimation of $f_1$ but that is significant. 

This time, the $p$-value of the Fisher's test $H_0:\beta_1^{(1)}=\hdots=\beta_d^{(1)}=0$ against $H_1:\exists j\in\{1,\hdots,d\}, \beta_j^{(1)}\neq 0$ is very small ($<2\times 10^{-4}$) which indicates that we are not close to a stationary point. Then, we try to improve the response doing a second descent step. 
  
\begin{figure}
\centering
\includegraphics[width=0.6\textwidth]{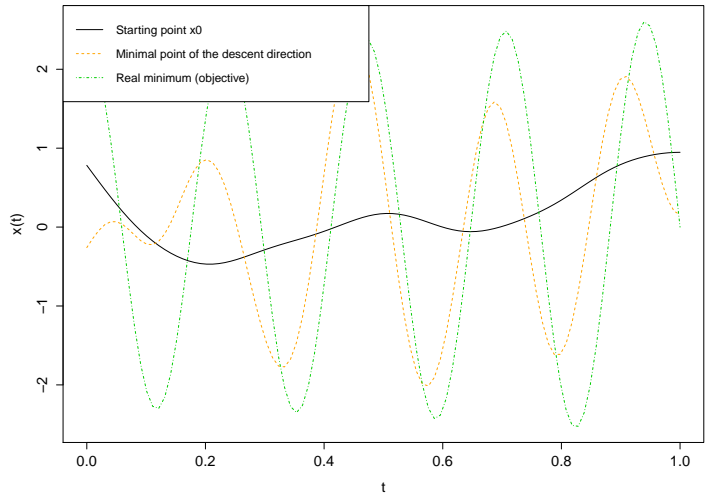}
\caption{\label{fig:result2} Result of optimization algorithm. Black curve: $x_0^{(0)}$, orange curve: $x_0^{(1)}$, red curve: $x_0^{(2)}$, green curve: $f_2$.}
\end{figure}

\subsection{Choice of basis}
\begin{figure}[h!]
\begin{tabular}{cc}
\includegraphics[width=0.45\textwidth]{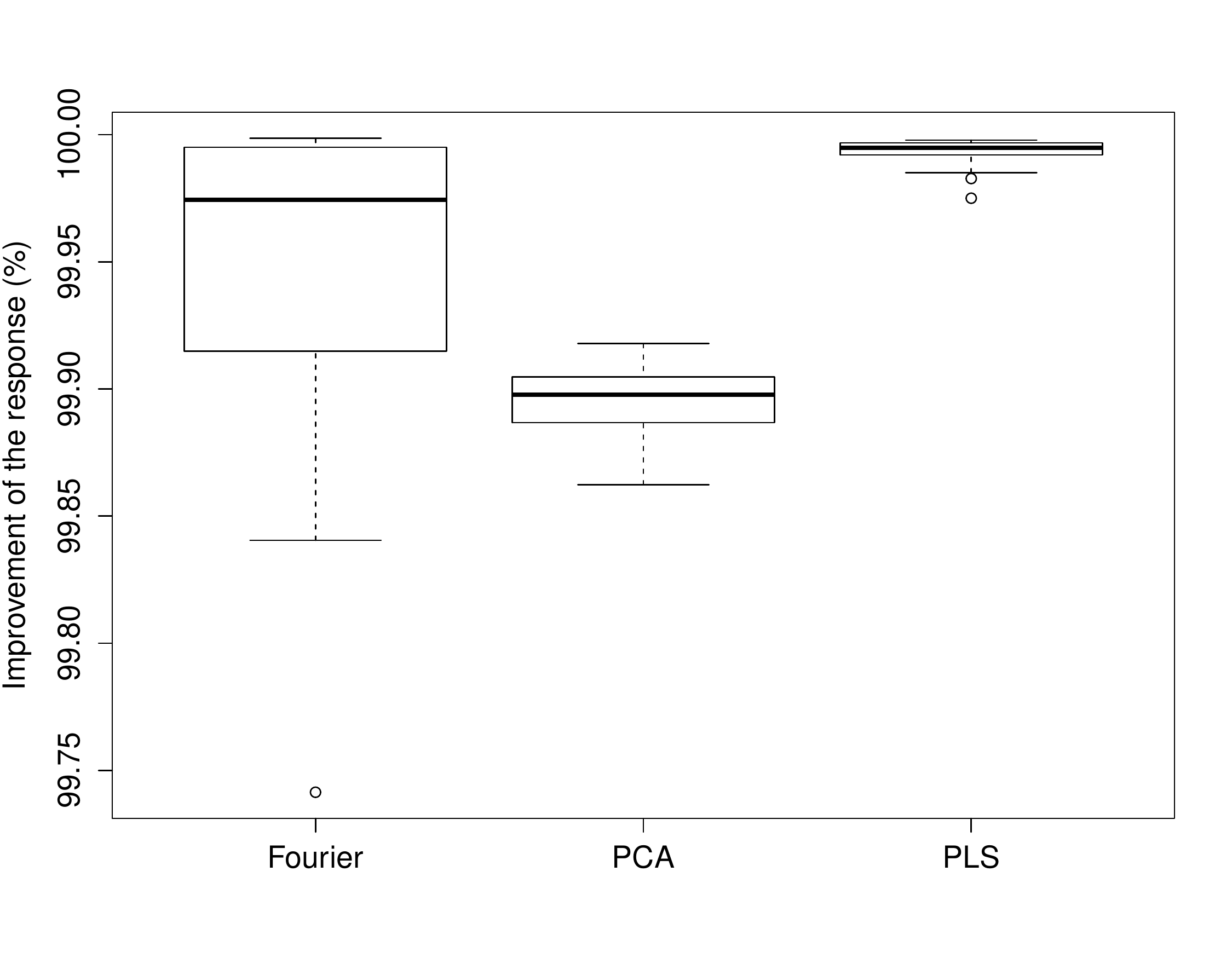}&\includegraphics[width=0.45\textwidth]{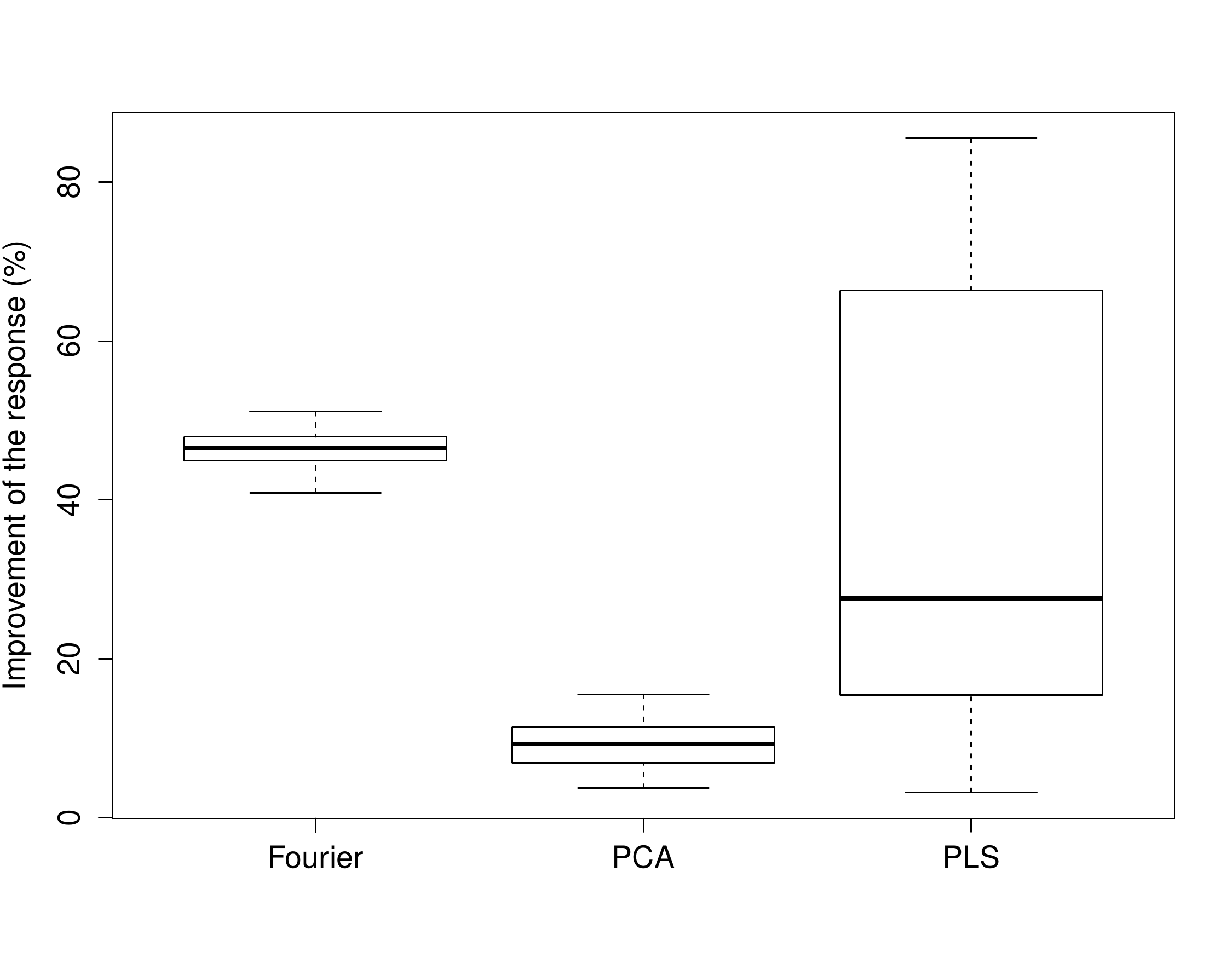}
\end{tabular}
\caption{\label{fig:MCstudy_bases}Monte-Carlo study of response improvement $\frac{m(x_0^{(0)})-m(x_0^{(1)})}{m(x_0^{(0)})}$ after the first descent step. Left-hand side: estimation of $f_1$, right-hand side: estimation of $f_2$.}
\end{figure}

In this section, we compare the three bases proposed in Section~\ref{sec:func_design}.
We generate $n_s=50$ training samples of size $n=500$ and compare the results of the first descent step when the design is generated by the Fourier basis, the PCA basis and the PLS basis. The starting point is the same: $x_0^{(0)}=X_{i_{\min}}$ for $i_{\min}={\arg\min}_{i=1,\hdots,n}\{Y_i\}$ (then for the Fourier basis the training sample is only used to set the starting point). The results are given in Figure~\ref{fig:MCstudy_bases}. We see immediately that the PLS basis seems to be a better choice than the PCA one. However, surprisingly, the choice between the PLS basis and the Fourier basis is less clear. 

\subsection{Choice of dimension $d$}

\begin{figure}
\begin{tabular}{c}
\includegraphics[width=\textwidth]{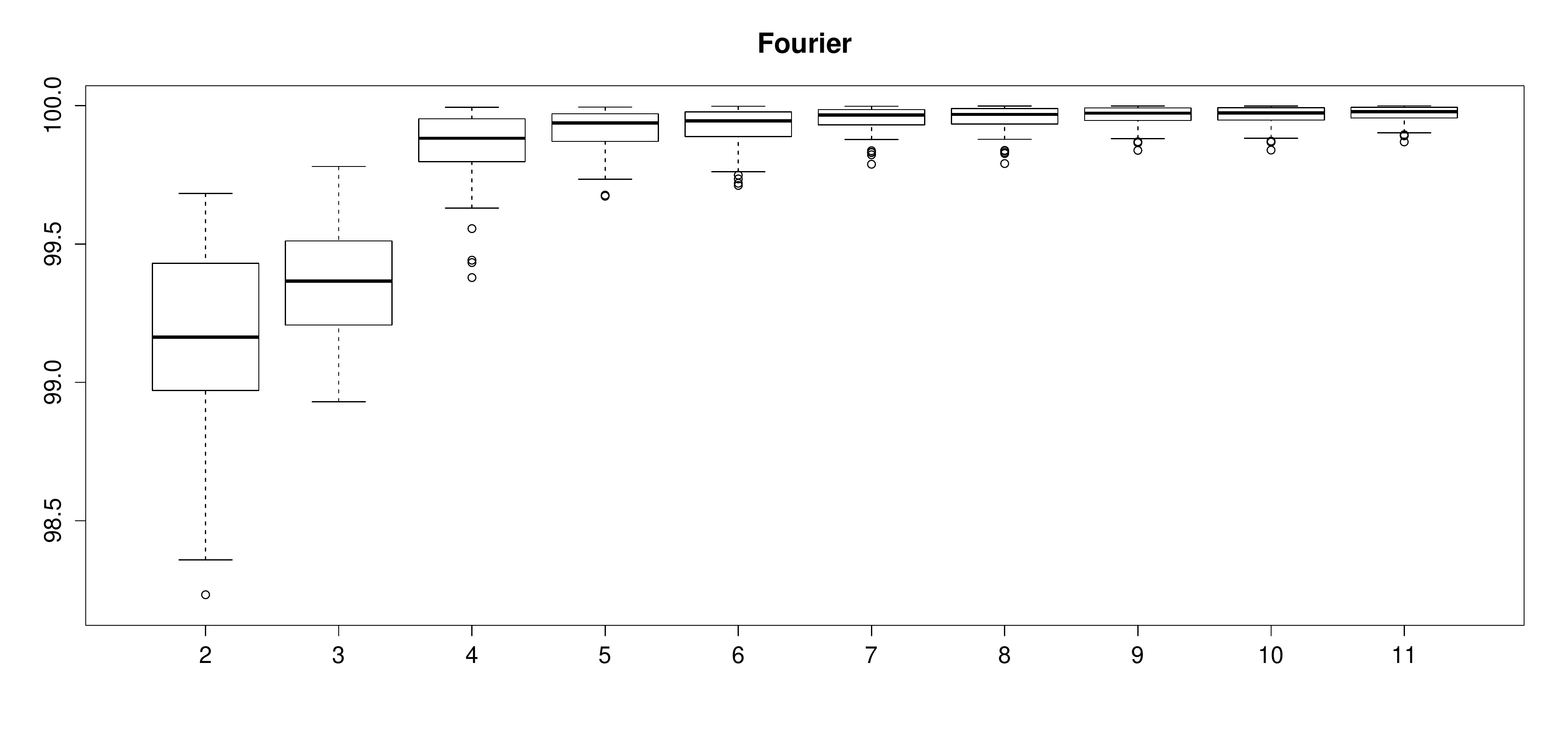}\\
\includegraphics[width=\textwidth]{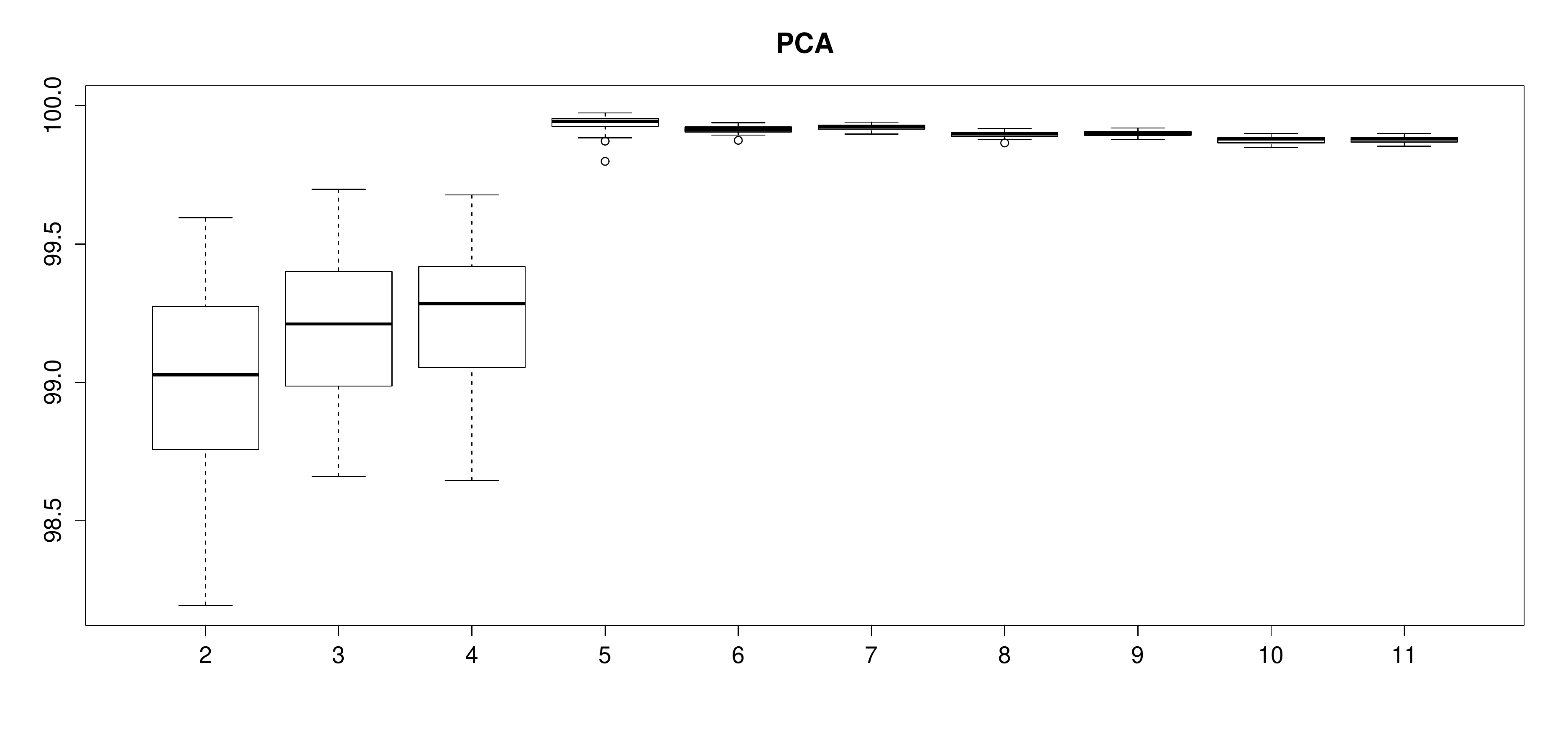}\\
\includegraphics[width=\textwidth]{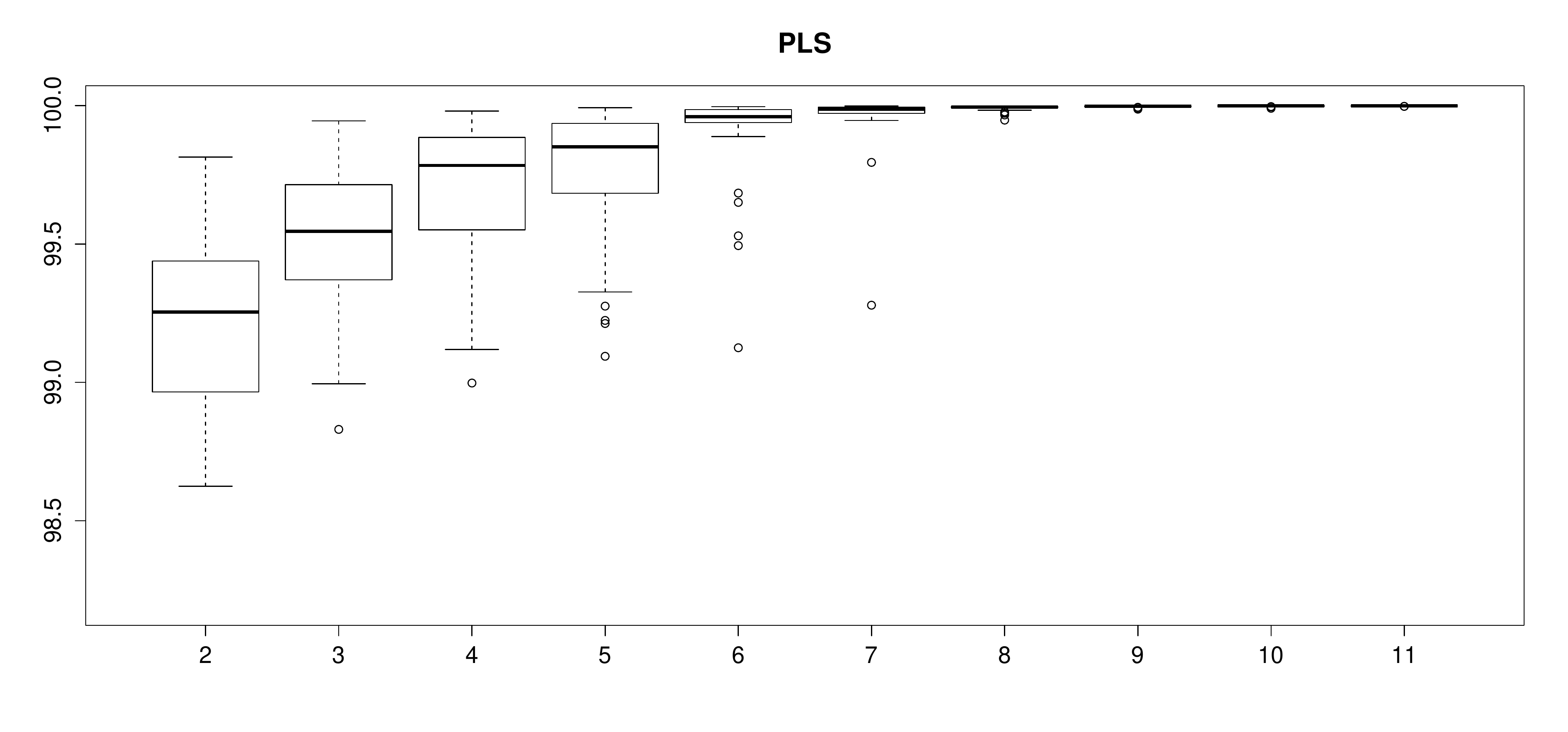}
\end{tabular}
\caption{\label{fig:MCstudy_dimension1}Monte-Carlo study of response improvement for the approximation of $f_1$ as a function of the dimension $d$. }
\end{figure}

\begin{figure}
\begin{tabular}{c}
\includegraphics[width=\textwidth]{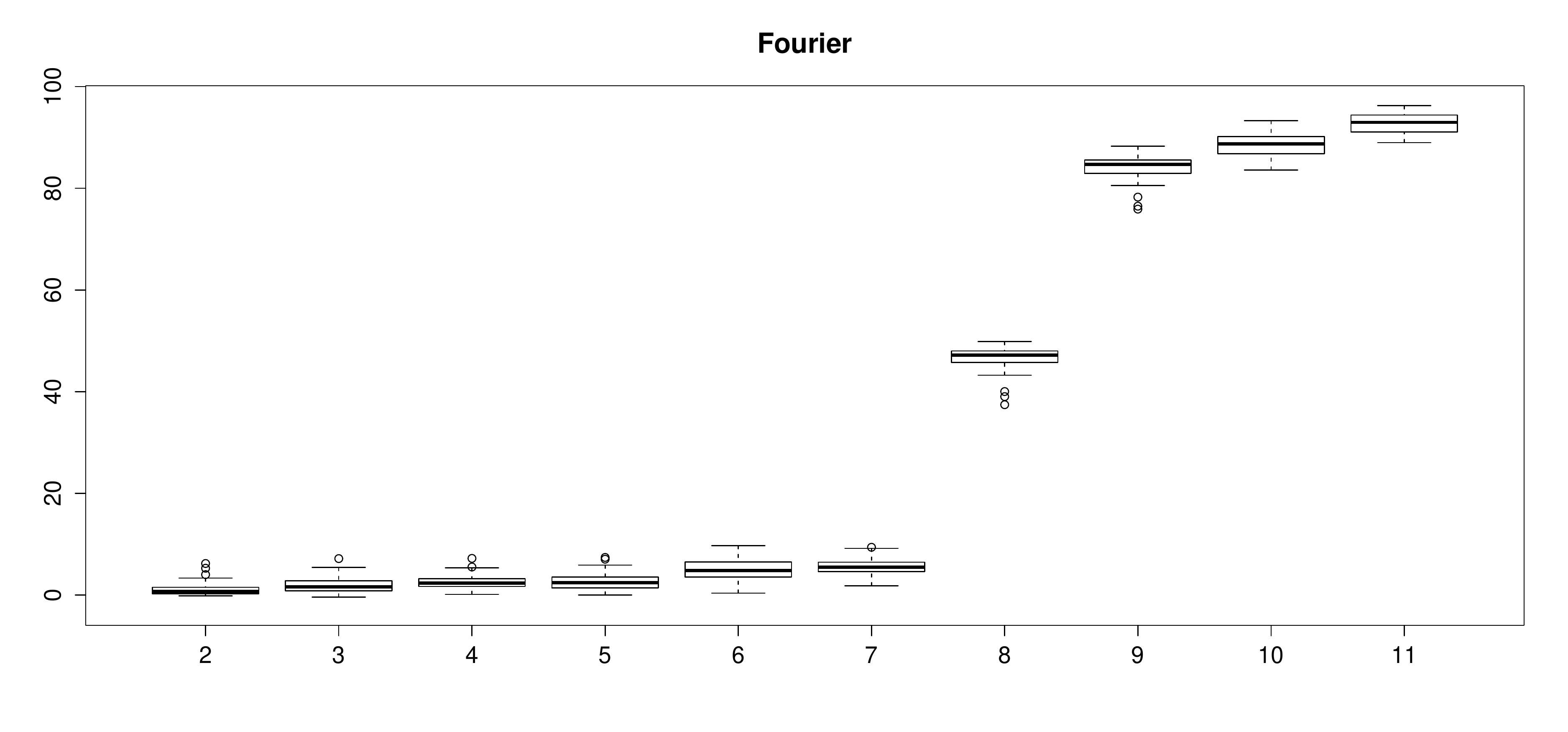}\\
\includegraphics[width=\textwidth]{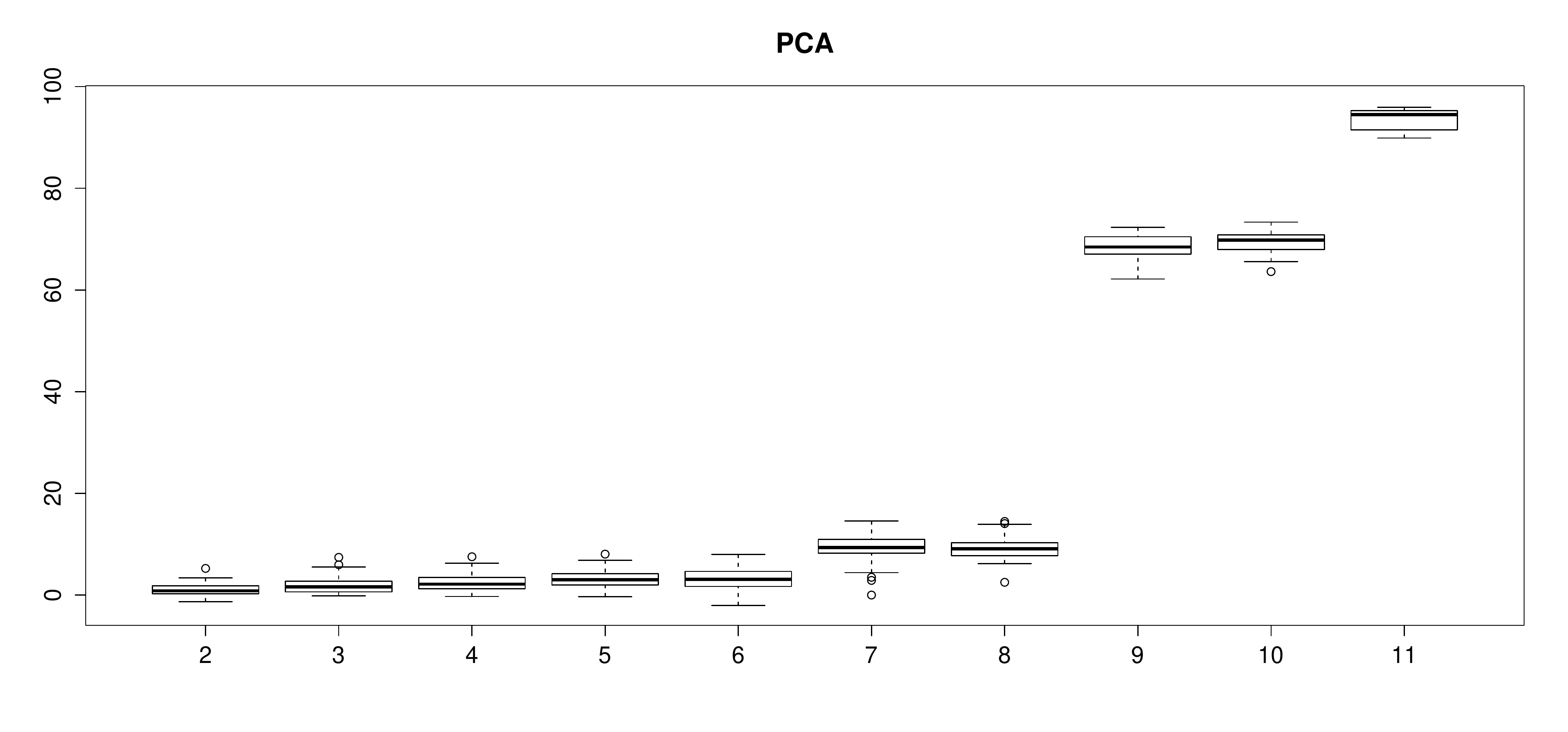}\\
\includegraphics[width=\textwidth]{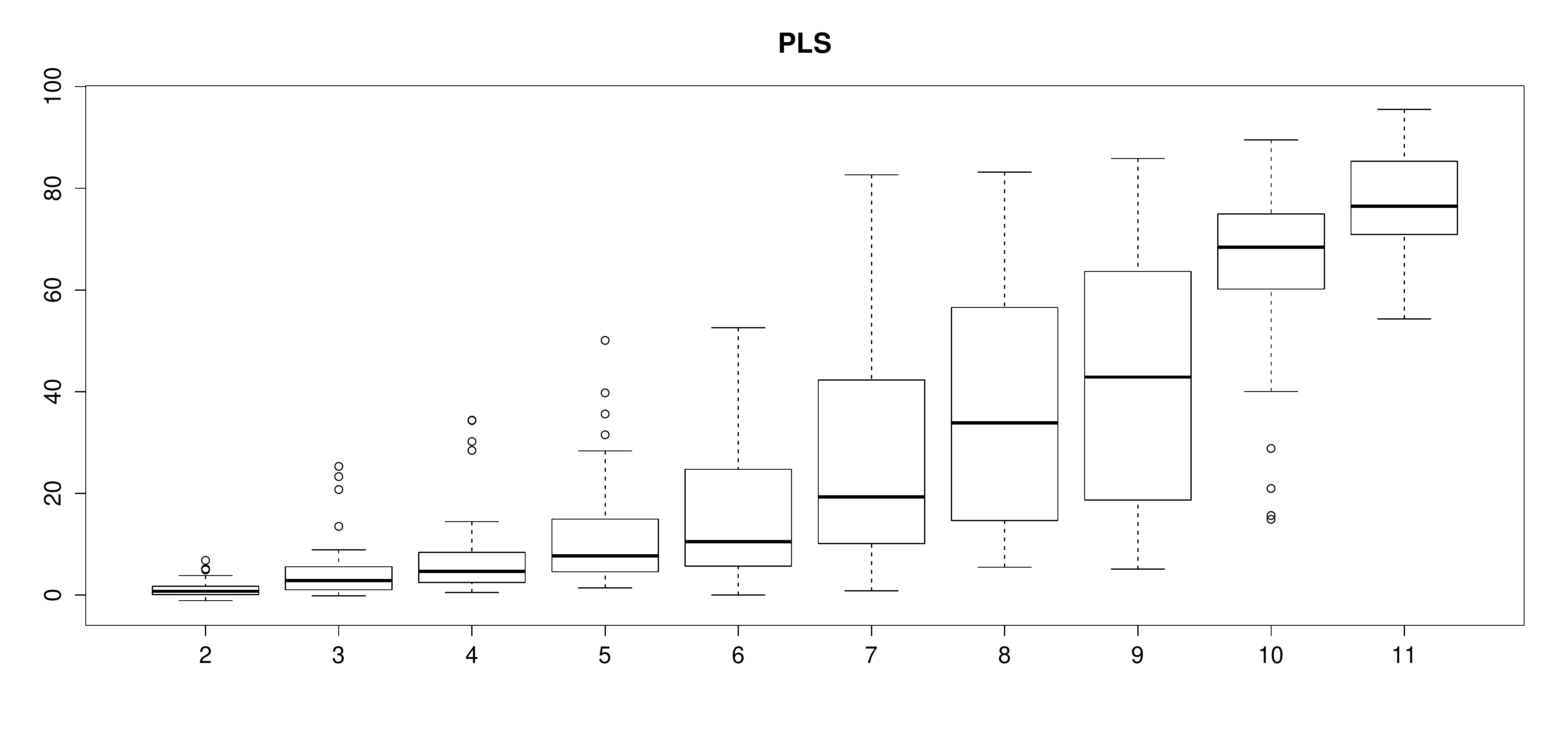}
\end{tabular}
\caption{\label{fig:MCstudy_dimension2}Monte-Carlo study of response improvement for the approximation of $f_2$ as a function of the dimension $d$. }
\end{figure}

Figures~\ref{fig:MCstudy_dimension1} and \ref{fig:MCstudy_dimension2} show that, except when the design is generated by the PCA basis for the approximation of $f_1$, the percentage of improvement increases when the dimension increases, which is coherent with the fact that the number of experiments grows exponentially with the dimension. 

\begin{figure}
\begin{center}
\includegraphics[width=0.5\textwidth]{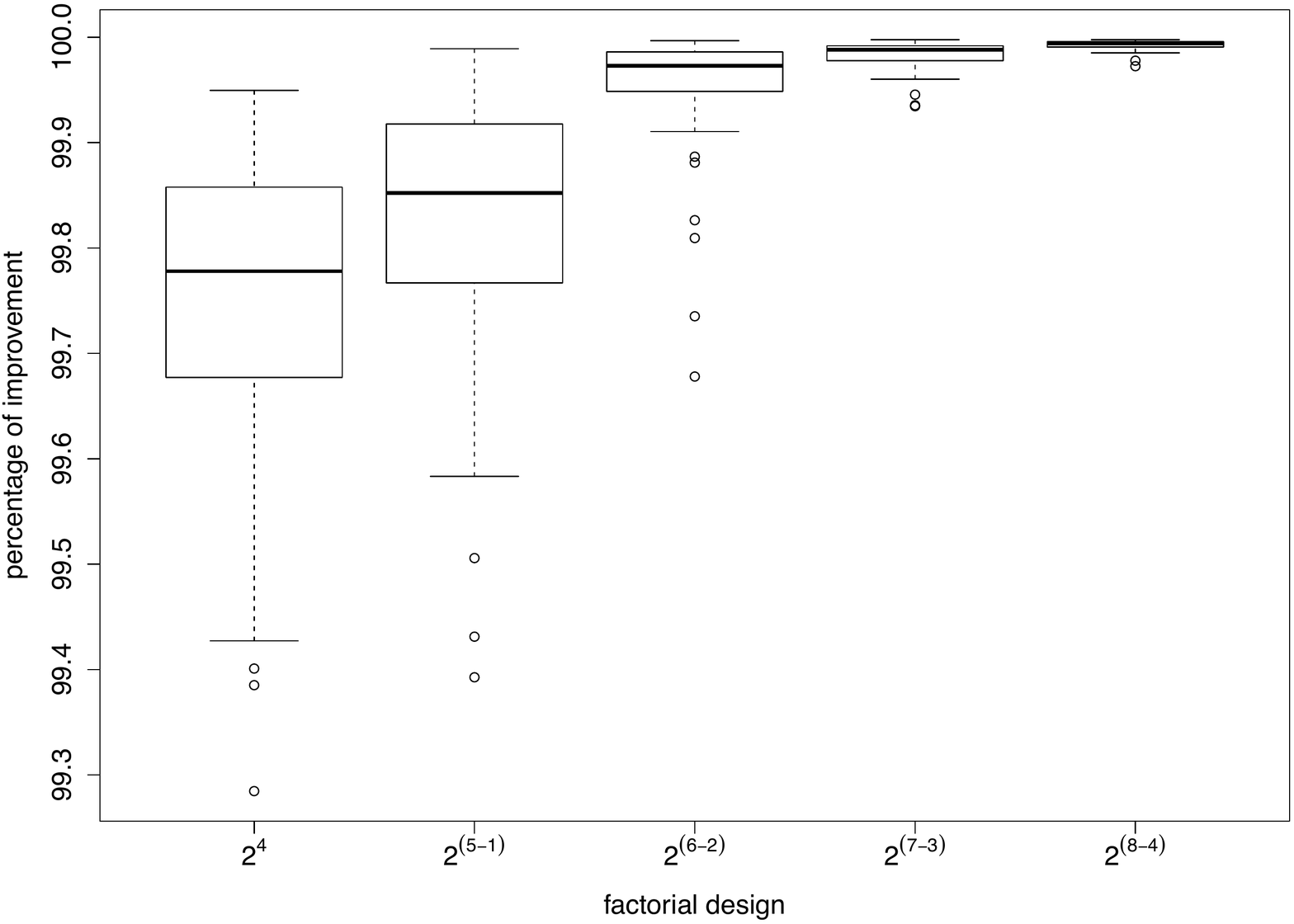}
\caption{Monte-Carlo study of response improvement for the approximation of $f_1$  with a $2^{d-p}$ factorial design, for different values of $d$ and $p$ such that $d-p=4$. }
\label{fig:etude_fact_design}
\end{center}
\end{figure}

We then decide to study the properties of the method when the number of design points is fixed and the dimension $d$ varies. For this purpose, we consider $2^{d-p}$ fractional factorial designs. We can see on Figure~\ref{fig:etude_fact_design} the results of the Monte-Carlo study. It seems that there is a significant improvement of the method when the dimension $d$ increases. This suggests that it is always better to explore new dimensions, even if we perform less experiments along each direction. 

\section{Application to nuclear safety}
\label{sec:CEA}

\subsection{Data and objectives}
An hypothetical cause of nuclear accident is the loss of coolant accident (LOCA). This is caused by a breach on the primary circuit. In order to avoid reactor meltdown, the safety procedure consists in  incorporating cold water in the primary circuit. This can cause a pressurised thermal shock on the nuclear vessel inner wall which increases the risk of failure of the vessel.

\begin{figure}
\begin{tabular}{ccc}
Temperature&Pressure&Heat transfer\\
\includegraphics[width=0.3\textwidth,height=4cm]{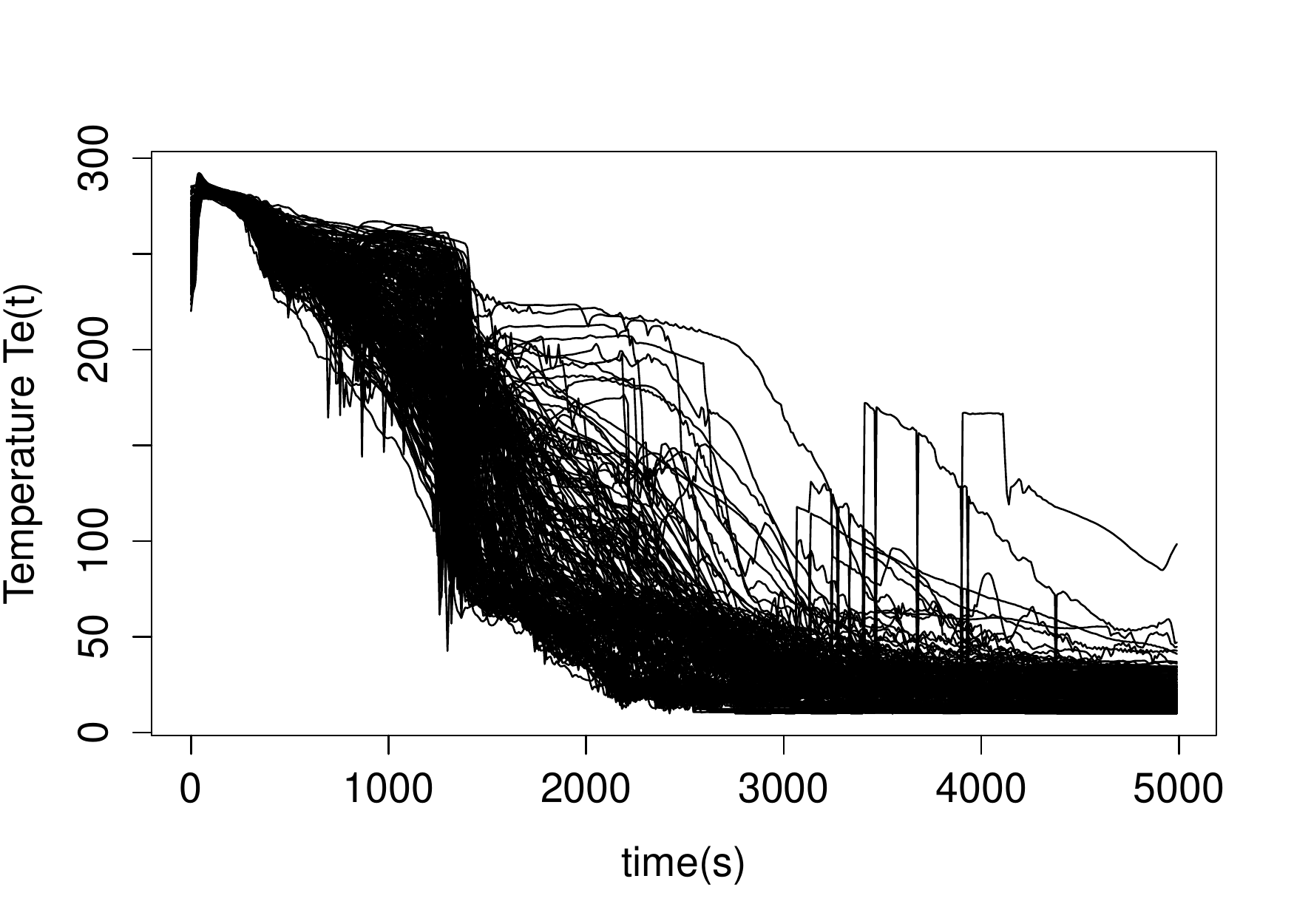}&\includegraphics[width=0.3\textwidth,height=4cm]{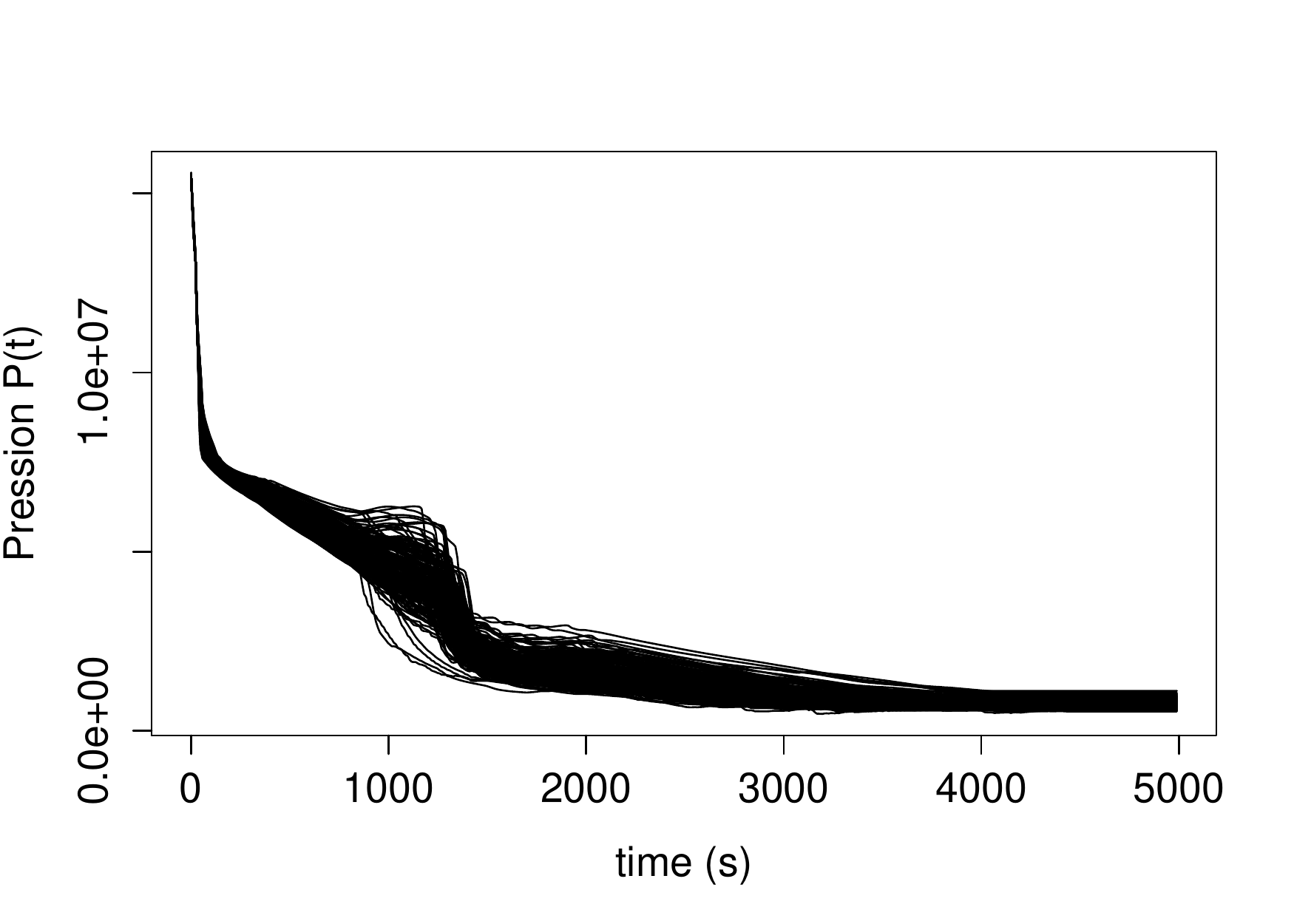}&\includegraphics[width=0.3\textwidth,height=4cm]{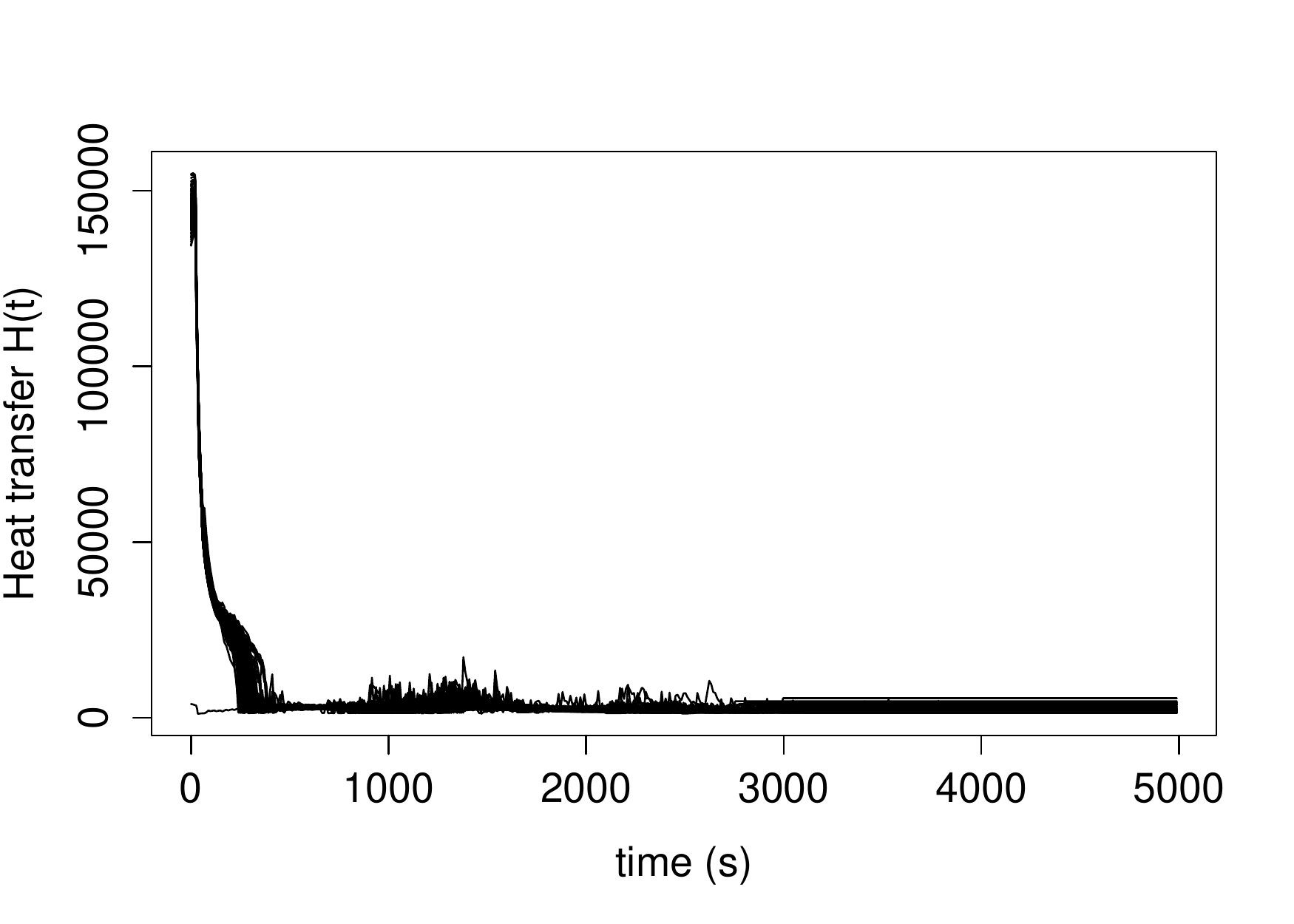}
\end{tabular}
\caption{\label{fig:courbes_transitoire_CEA}Evolution of temperature, pressure and heat transfer (learning sample). Source: CEA.}
\end{figure}

The parameters influencing the probability of failure are the evolution over time of temperature, pressure and heat transfer in the vessel. Obviously, the behavior of the reactor vessel during the accident can be hardly explored by physical experimentation and numerical codes have been developed, for instance by the  CEA\footnote{French Alternative Energies and Atomic Energy Commission (Commissariat \`a  l'\'energie atomique et aux \'energies alternatives), government-funded technological research organisation.  \url{http://www.cea.fr/}}, reproducing the mechanical behavior of the vessel given the three mentioned parameters (temperature, pressure, heat transfer). Figure~\ref{fig:courbes_transitoire_CEA} represents different evolution of each parameter during the procedure depending on the value of several input parameters, which can be used as a learning sample.

The aim is to find the temperature transient which minimizes the risk of failure. We have access here to the margin factor (MF) which decreases when the risk of failure increases. Hence, the aim is to maximise the MF.  
  
\subsection{Generation of design}

Considering that the inputs are the three temperature, pressure and heat penetration curves, we set $\mathbb H=\left(\mathbb L^2([0,T])\right)^3$ (with $T=5000 s$) equipped with the natural scalar product 
\begin{eqnarray*}
\langle (x^{(T)}_1,x^{(P)}_1,x^{(H)}_1),(x^{(T)}_2,x^{(P)}_2,x^{(H)}_2)\rangle&=&\\
&&\hspace{-4cm}\int_0^T x^{(T)}_1(t)x^{(T)}_2(t)dt+\int_0^T x^{(P)}_1(t)x^{(P)}_2(t)dt+\int_0^T x^{(H)}_1(t)x^{(H)}_2(t)dt.
\end{eqnarray*}

We define the starting point of the algorithm as the triplet $(X_i^{(T)},X_i^{(P)},X_i^{(H)})$ of the learning sample maximizing the response. 

\begin{figure}
\begin{tabular}{ccc}
(a)&(b)&(c)\\
\includegraphics[width=0.3\textwidth]{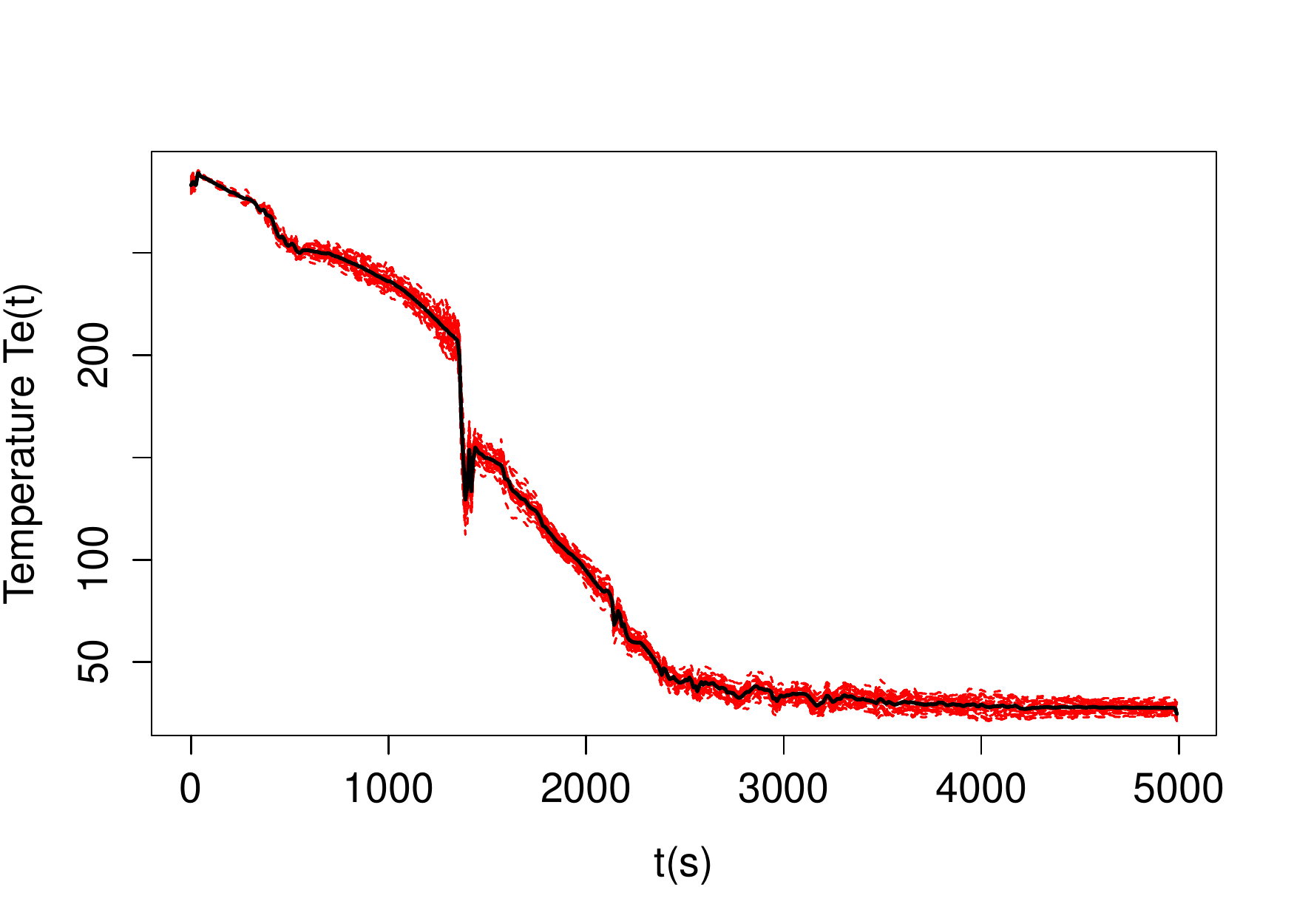} &\includegraphics[width=0.3\textwidth]{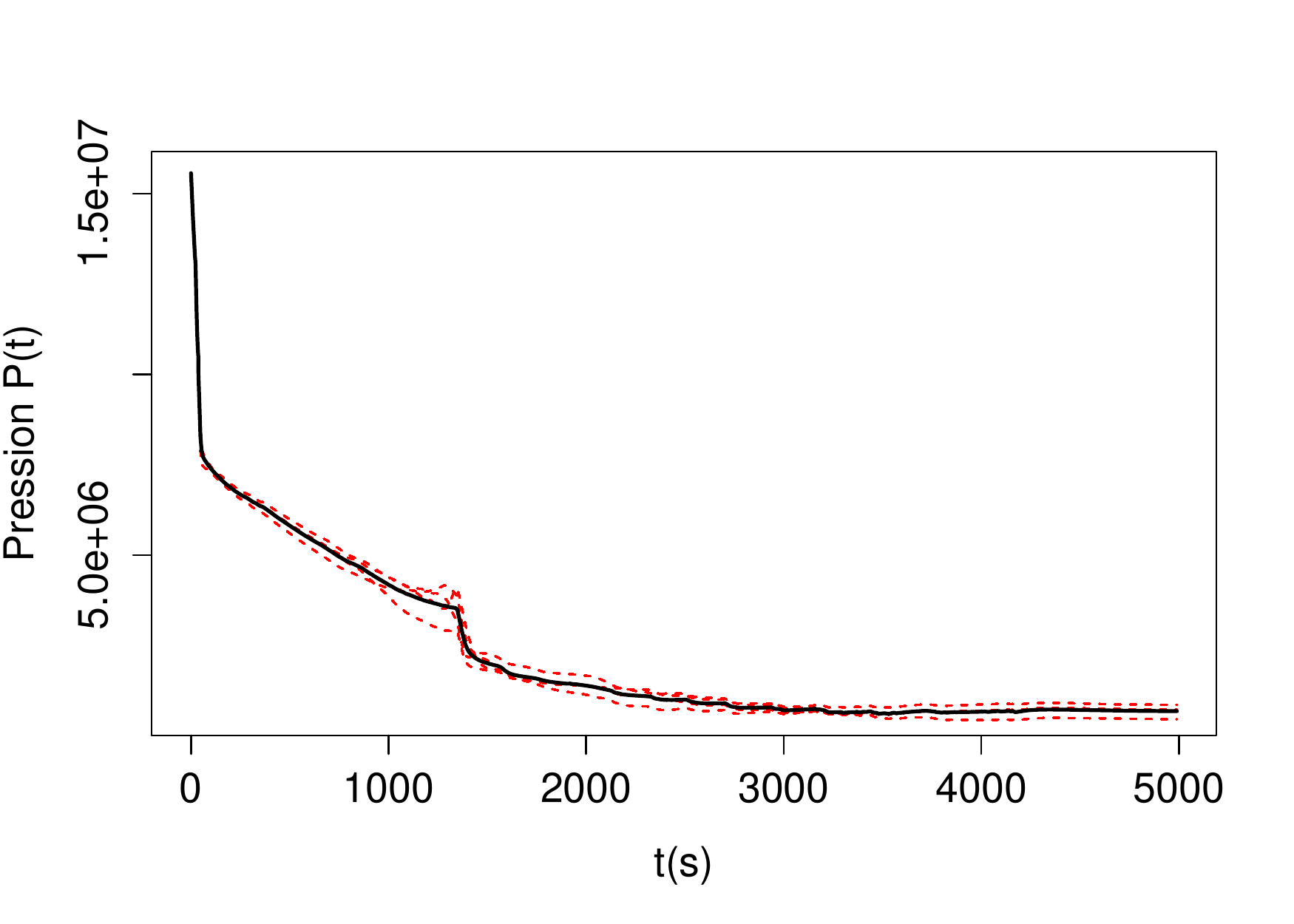} &\includegraphics[width=0.3\textwidth]{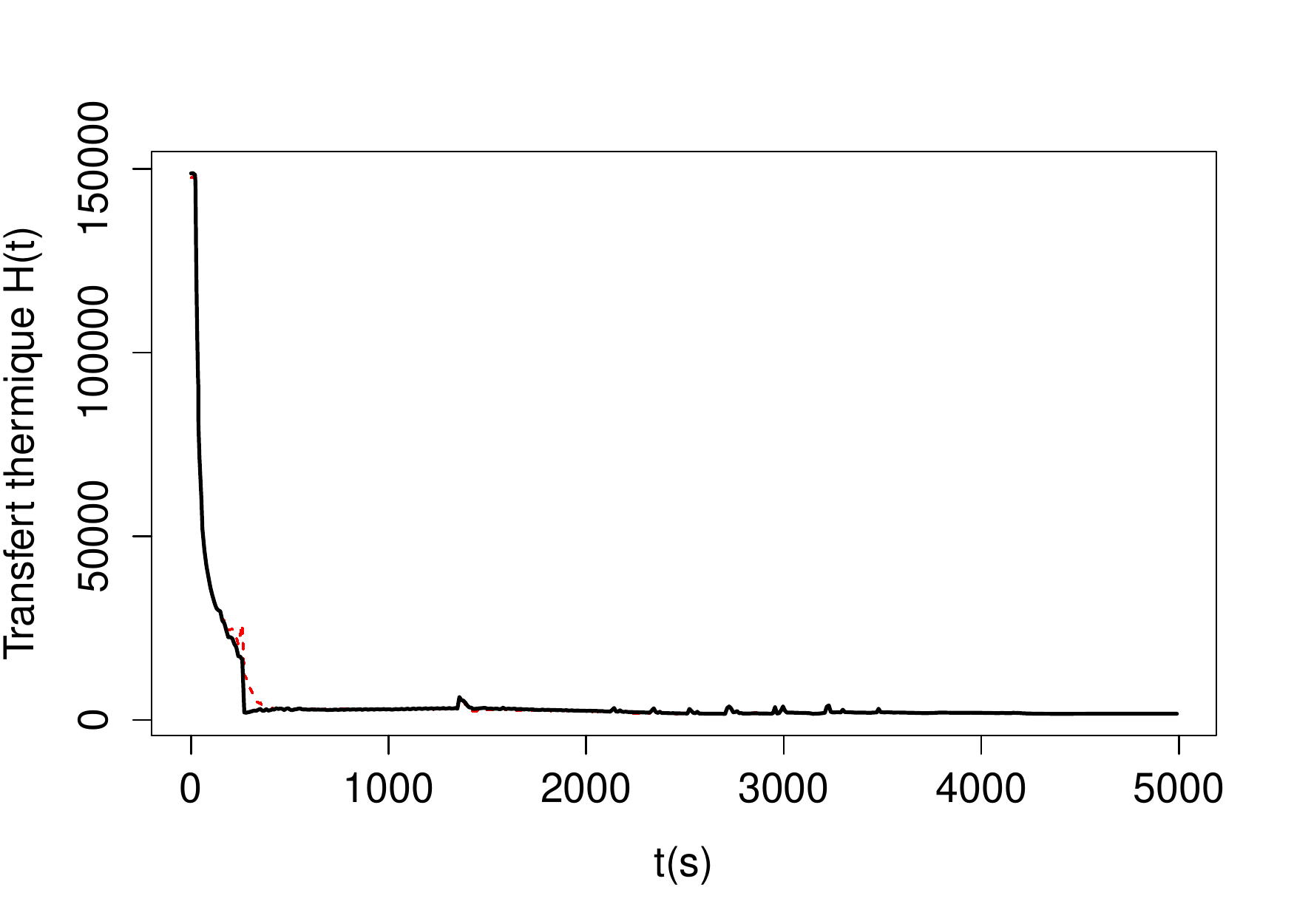} \\
\end{tabular}
\caption{\label{fig:design_CEA}Functional experimental design around the initial curves.  }
\end{figure}

In view of the simulation results of Section~\ref{sec:num_exp} and the presence of a learning sample, we focus on the PLS basis and generate a functional design based on a minimum aberration $2^{10-5}$ fractional design for the temperature, a $2^{3-2}$ design for the pressure and the heat transfer.  As some design points of the functional design around the initial heat transfer curve took negative values (which can not correspond to the physic since the heat transfer is always positive), we remove it and keep only the design points which are always positive. The design points are plotted in Figure~\ref{fig:design_CEA}. The resulting design, which is a combination of all curves of the three designs obtained (for temperature, pressure and heat penetration) counts 128 design points.   

\subsection{Results}

\begin{figure}
\begin{tabular}{cc}
\includegraphics[width=0.45\textwidth]{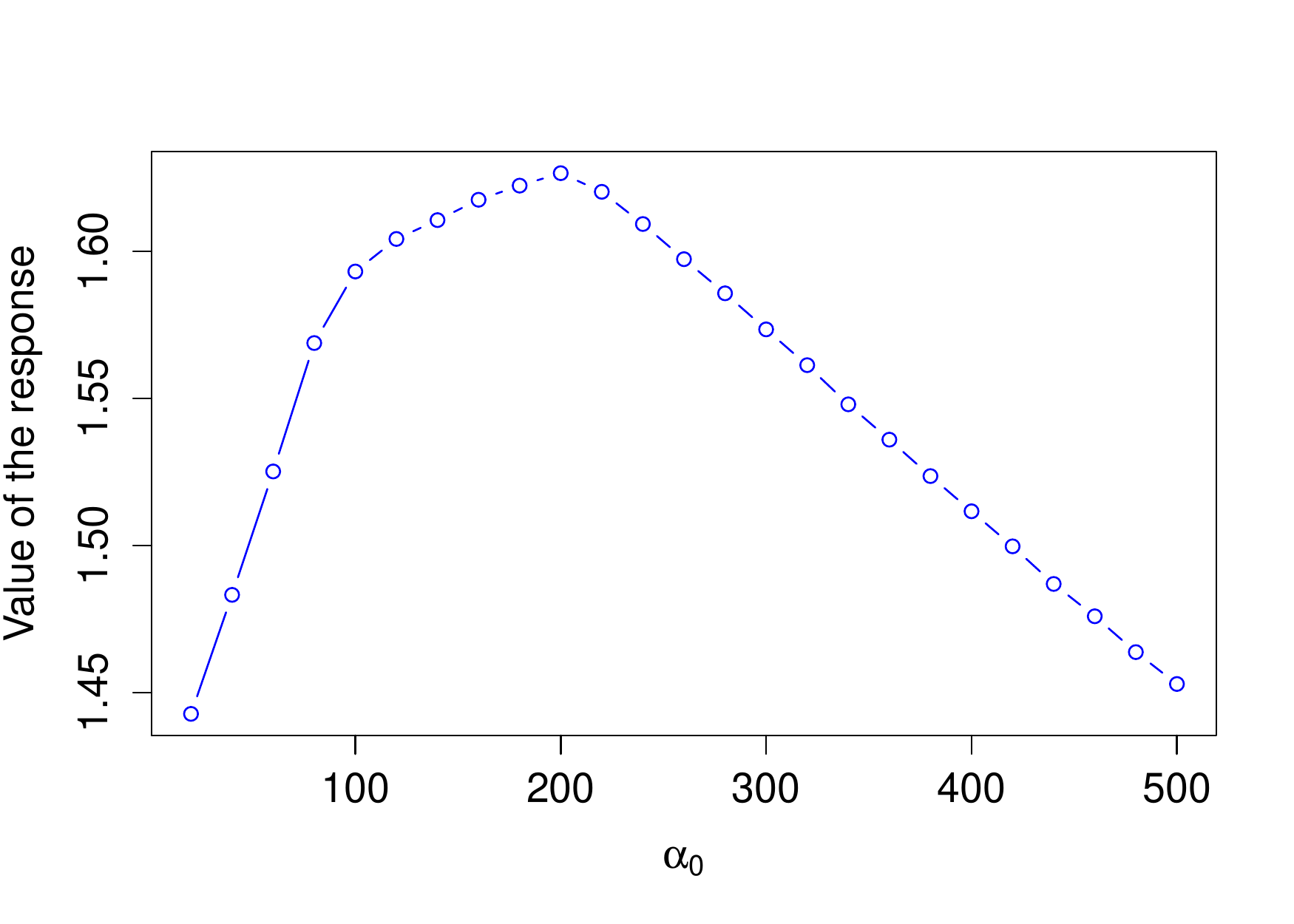}&\includegraphics[width=0.45\textwidth]{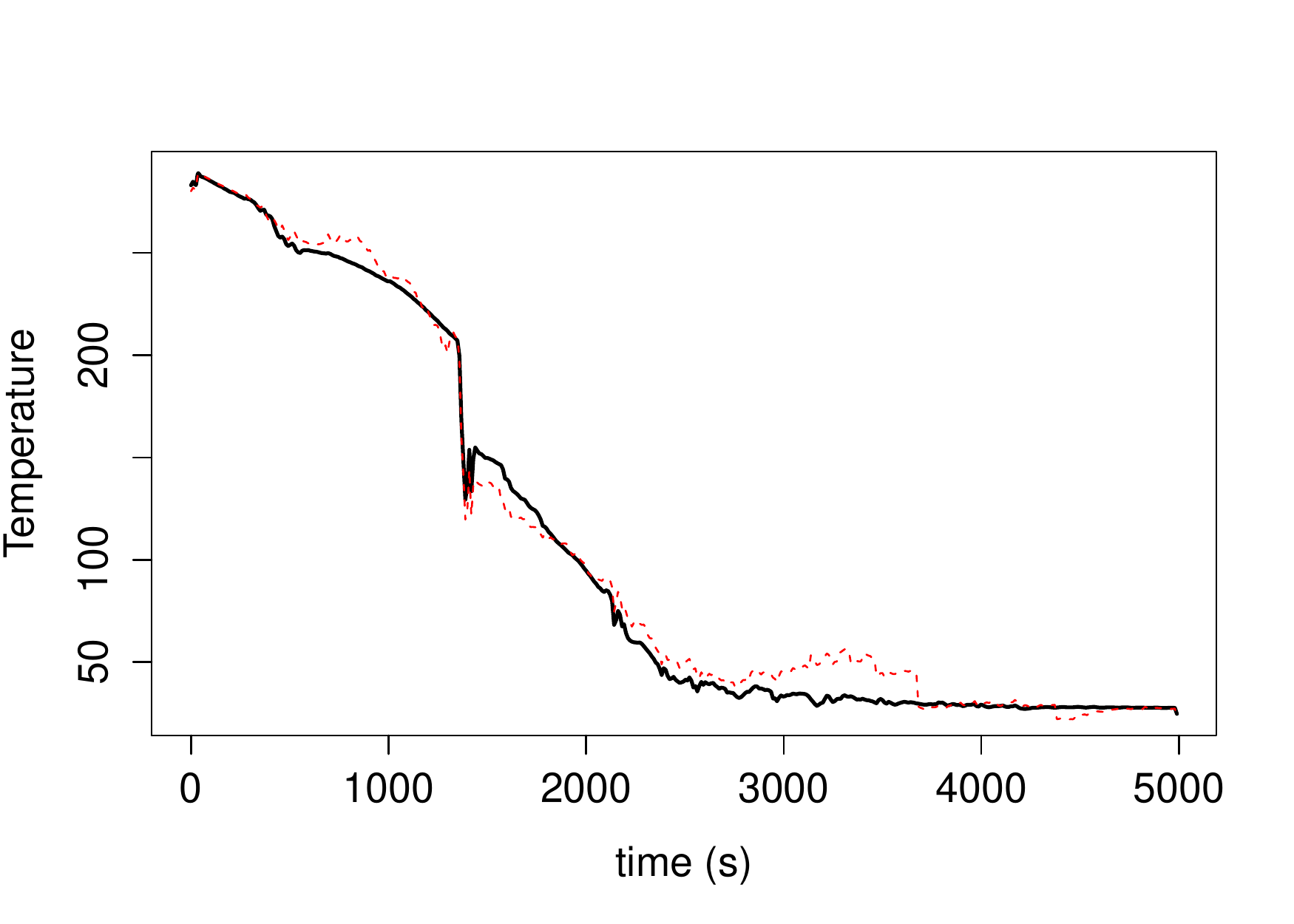}
\end{tabular}

\caption{\label{fig:dire_steep_CEA}Left: value of the response on the estimated steepest ascent direction. Right: solid line initial temperature point, dotted line: optimal temperature transient estimated. }
\end{figure}

\begin{figure}
\begin{tabular}{ccc}
\includegraphics[width=0.3\textwidth]{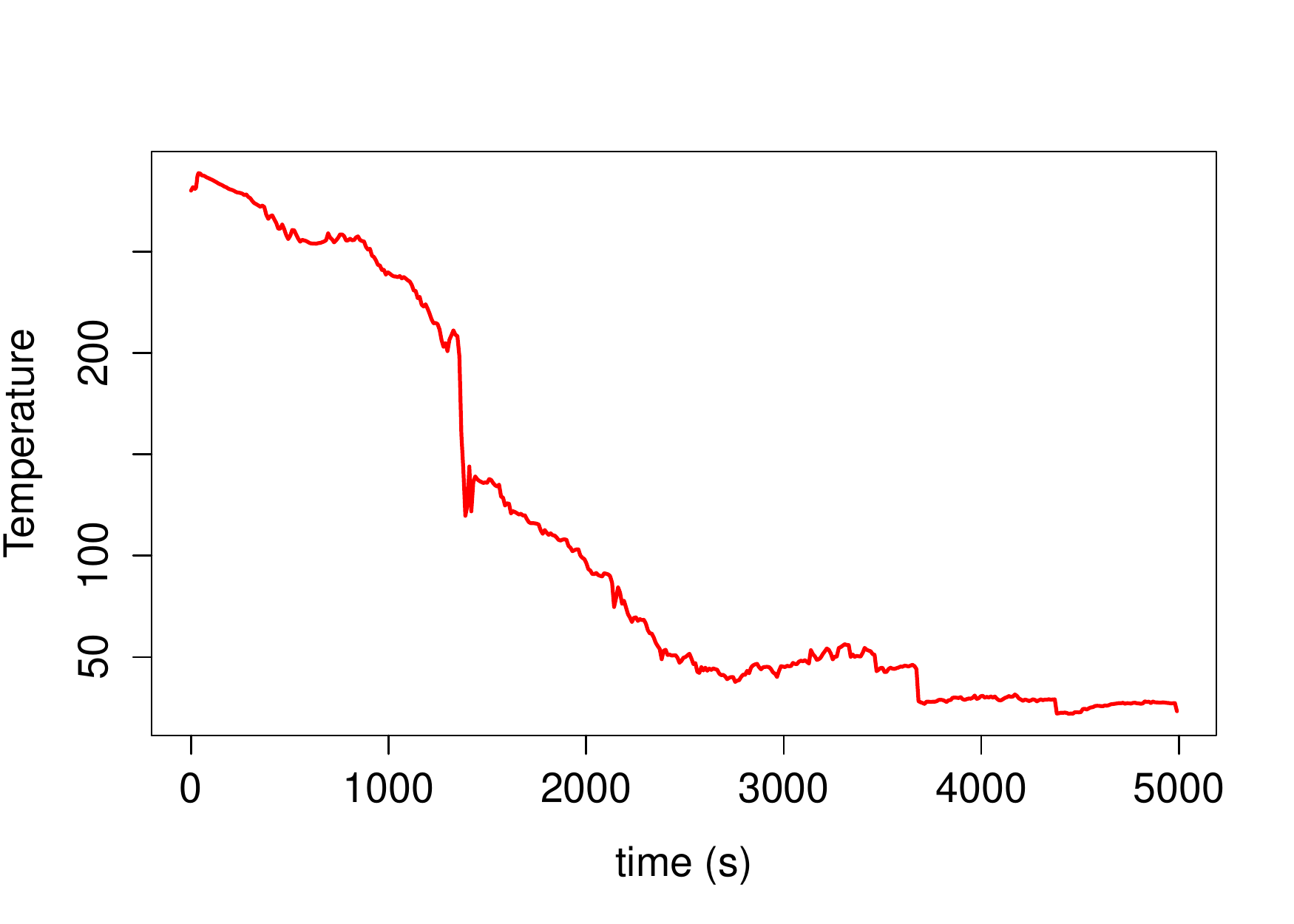}&\includegraphics[width=0.3\textwidth]{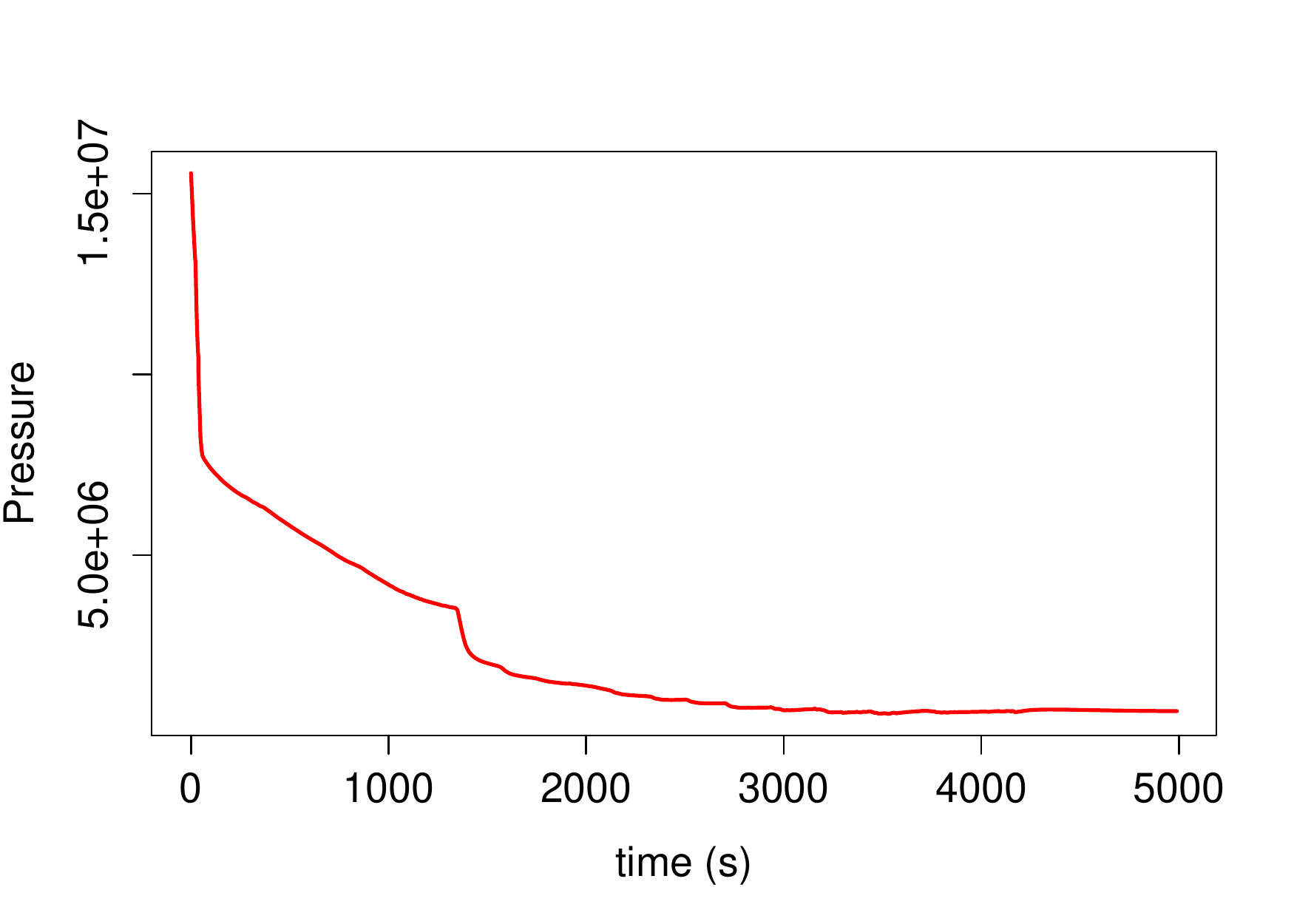}&\includegraphics[width=0.3\textwidth]{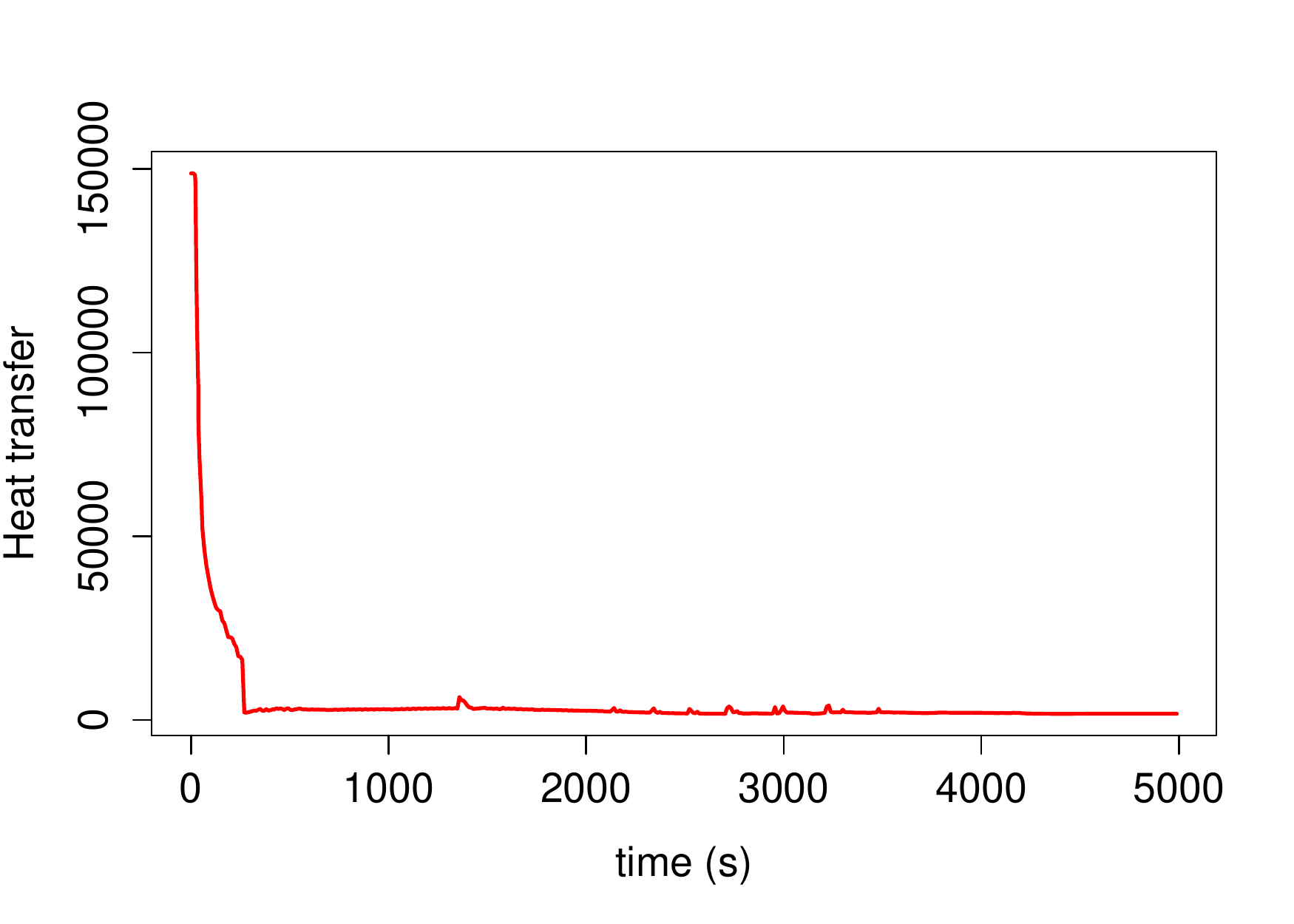}
\end{tabular}
\caption{\label{fig:final} Point of the estimated steepest ascent direction maximizing the response.}
\end{figure}

We compute an estimation of the gradient with the results of the experiments on the design points given in Figure~\ref{fig:design_CEA}. The results are given in Figure~\ref{fig:dire_steep_CEA}. We take $\lambda_0=200$. The final estimates of the optimal curves are given in Figure~\ref{fig:final}.  

The main change between the starting point and the minimal point of the ascent direction lies in the temperature transient (the changes in the pressure and heat temperature transient are minimal), especially in the evolution of temperature between 500\,s and 2000\,s after the simulation starts and between 3000\,s and 4000\,s.  

\section*{Acknowledgement}
The data were provided by the French Alternative Energies and Atomic Energy Commission (CEA). The author wants to thank especially Michel Marques for its patience and cooperation, which was essential to obtain the final results. I also want to thank \'Elodie Brunel and Andr\'e Mas for their helpful advices and careful reading of this work as well as Herv\'e Cardot and Peggy C\'enac for their precious help on stochastic optimization. 
\bibliographystyle{abbrvnat}
\bibliography{biblio_these} 
\end{document}